\documentclass{amsart}
\usepackage{times}
\usepackage{amsthm}
\usepackage{amssymb}
\usepackage[pdftex]{graphicx}
\usepackage[all]{xy}
\usepackage[mathscr]{eucal}
\usepackage{amsmath,latexsym,oldlfont}
\usepackage{mathdots}
\usepackage{pifont}
\usepackage{stmaryrd}
\usepackage{textcomp}
\usepackage{pifont}
\usepackage{multirow}
\usepackage[colorlinks,linkcolor=blue,citecolor=red,linktocpage=true]{hyperref}
\usepackage{mathtools}
\numberwithin{equation}{section}
\newtheorem{teorema}{Teorema}[section]

\newtheorem{theorem}[teorema]{Theorem}
\newtheorem{lemma}[teorema]{Lemma}

\newtheorem{remark}[teorema]{Remark}
\newtheorem{corollary}[teorema]{Corollary}
\newtheorem{definition}[teorema]{Definition}

\newsavebox{\spacebox}
\begin{lrbox}{\spacebox}
    \verb*! !
\end{lrbox}

\title[Iterative methods to build LG-matrices and  Applications]{Iterative methods to build LG-matrices and  Applications}
\author[J. Carrillo--Pacheco]{Jes\'us Carrillo--Pacheco}
\author[F. Jarqu\'{\i}n--Z\'arate]{Fausto Jarqu\'{\i}n--Z\'arate}
\address[J. Carrillo--Pacheco and F. Jarqu\'{\i}n--Z\'arate]{Academia de Matem\'aticas, Universidad Aut\'onoma de la Ciudad de M\'exico, 09390 Ciudad de M\'exico,  M\'exico.}
\email[J. Carrillo--Pacheco]{jesus.carrillo@uacm.edu.mx}
 \email[F. Jarqu\'{\i}n--Z\'arate]{fausto.jarquin@uacm.edu.mx}
 \thanks{Jes\'us Carrillo-Pacheco and Fausto Jarqu\'{\i}n-Z\'arate  were supported by Laboratorio de
Cifrado y Codificaci\'on de la Informaci\'on (LCCI-UACM) of
Universidad Aut\'onoma de la Ciudad de M\'exico.}

\begin{document}



\begin{abstract}
{In this paper  we give a recursive algorithm to construct two families of $(0,1)$-matrices, 
one sparse regular and the other dense.
We study various properties of the two families of $(0,1)$-matrices built with our algorithm.
We present a new construction of two clases of isodual linear codes, 
one is the low density generator matrix codes and other is the dense linear codes, 
for both codes we obtain the polynomial of the distribution of weights, a bound for the minimum distance and
we apply to these codes the efficient encoders based on approximate lower triangulations developed 
by Richardson-Urbanke.
We identify the unique $(0,1)$-matrices, up  basis change, associated with the geometry 
of the Lagrangian-Grassmannian variety.}
\end{abstract}

 \maketitle

\section{Introduction}\label{Secc0}

 Matrices with special ``structure''  arise in many parts of mathematics and its applications.
 Inspired by the geometry of the Lagrangian-Grassmannian variety we studied a type of LG-matrices.
More explicitly, the Lagrangian-Grassmannian variety $L(n,2n)$ parametrizes the set of all maximal isotropic subspaces 
(Lagrangian) of a symplectic vector space of dimension $2n$ over arbitrary field $\mathbb{F}$, 
this variety is a linear section of the Grassmann variety and we can see as a set of zeros of linear  homogeneous polynomials (see [3]) and has a type of LG-matrix associated. 
Following this path the reader can see a connection between various areas of mathematics, the symplectic geometry,
theory of matrices ( (0,1)-matrix, persymmetric matrix, incidence matrix) and applications.

In this paper, we present a recursive algorithm to construct two families of $(0,1)$-matrices.
The first family of sparse regular  type,  for this,  we start  with $k$ and $\ell$ any positive integers and we build a structured sparse $(0,1)$-matrix  $A_k^{\ell-2}$ in approximate lower triangular form 
of order $C_{\ell-1}^{k+\ell-1} \times C_\ell^{k+\ell-1}$,
with $k$ ones in each row and $\ell$ ones in each column for $k$ and $\ell$ any positive integers.

 The second family of dense type, that is, modifying the algorithm of the first family we construct $B_k^{\ell-2}$ dense
 matrix of order   $C_{\ell-1}^{k+\ell-1} \times C_\ell^{k+\ell-1}$  with the property 
 $A_k^{\ell-2} + B_k^{\ell-2} = \mathbb{J}$, where $\mathbb{J}$ is the matrix full of ones. 
 The importance of this is that we propose an algorithmic construction.

We proof that
\begin{enumerate}
\item[i)] $A_k^{\ell-2}$ is  a sparce matrix and is in approximate lower diagonal form, 
$B_k^{\ell-2}$ is a densa matrix and $A_k^{\ell -2} + B_k^{\ell-2} = \mathbb{J}$, where $\mathbb{J}$ is the matrix full of ones.
\item[ii)]  $A_k^{\ell-2}$ and $B_k^{\ell-2}$ are fragmented matrices.

\item[iii)] $A_k^{k-2}$ and $B_k^{k-2}$ are persymmetric matrices.

\item[iv)] $A_k^{k-3}$ is an incidence matrix between a triangulated set $T_m$ and the set $\{0,1\}^{I(m/2,m)}$, 
where $k=\frac{m+2}{2}$ and $m$ is even positive integer.

\item[v)] Using $A_k^{k-3}$, we obtain $\mathcal{M}$ the incidence structure of section lineal of LG variety.
\end{enumerate}

This paper contains the following applications:
\begin{enumerate}
\item We present  two clases of isodual linear codes, 
one is the low density generator matrix codes and other is the dense linear codes, 
for both isodual codes we obtain the polynomial of the distribution of weights and  a bound for its the minimum distance. 
An important fact is that they have the properties of persymmetric and  similarity.

\item We introduce a family of Low Density Parity Check codes (LDPC-codes), 
these are a generalization of the codes considered in [4]. For our new family of codes we propose
an optimal decoding algorithm using ``Efficient encoders based on approximate lower triangulations'' developed by 
Richardson-Urbanke see [7]. Our construction has enough reduction of computation time to implement the method of
Richardson-Urbanke, this is due to  the matrix associate to LDPC-codes is in approximate lower triangular form.

\item We identify the unique $ (0,1) $-matrices $\mathcal{M}$, up  basis change, associated with the geometry of the Lagrangian-Grassmannian variety as a linear section of the Grassmann variety.
For this purpose, the combinatorial study of the $A_k^{k-3}$ matrix is essential, this is a technical part of the work, and
we proof that $\mathcal{M}$ is a block diagonal matrix. In the case that the scalar field $\mathbb{F}$ is finite, 
$\mathcal{M}$ uniquely characterizes the linear section of the Grassmann variety such that annuls the set of its rational 
points.
\end{enumerate}

The paper is organized as follows. Section \ref{Secc1}  contains preliminary and basic $(0,1)$-matrix notation. 
Section \ref{Secc2} and \ref{Secc3} we present two algorithm to construct the matrices $A_k^{\ell-2}$ and $B_k^{\ell-2}$ and we proof 
several of the properties that these have. Section \ref{Secc4} we shown that the matrix $A_k^{k-3}$ is incidence matrix of the configuration of triangle-sets $T_m$ of the set of indexes. Section \ref{Secc5} we proof that a sparse $(0,1)$-matrix 
$\mathcal{M}$ it is expressed as a direct sum of matrices $A_k^{k-3}$ see Theorems \ref{thmPeven} and \ref{thmPodd}. Section \ref{Secc6} contains the applications (a), (b) and (c). 
Finally the paper contains an Appendix that deals with the Lagrangian-Grassmannian variety.

\section{Preliminary}\label{Secc1}
If  a matrix $A$ has all its coefficients equal  $0$ or $1$  is called a $(0,1)$-{\it matrix}. Give a $(0,1)$-matrix $A$ we say that  is {\it regular} if the number of 1's is fixed in each column and 
has a fixed number of 1's in each row. If $A$ is not regular we say that  is {\it irregular}. 
We say that $H$ is a matrix {\it approximate lower diagonal form}  if in the lower left corner of  $H$ is has an identity matrix as its submatrix.
Following  \cite{5}, let $A$ be a $(0, 1)$-matrix of order $p \times q$,  the sum of row $i$ of $A$ be denoted by $r_i$ where  
$i=1,\dots, p$ and  the sum of column $j$ of $A$ be denoted by $s_j$ where  $j = 1, \ldots, q$. 
With the $(0,1)$-matrix $A$ we associate the row sum vector $R=(r_1,\ldots, r_ p)$
where the $i$-th component gives the sum of row $i$ of $A$. Similarly, the column sum vector $S$ is denoted by
$S=(s_1, \ldots, s_q)$. The vectors $R$ and $S$ determine a class ${\EuScript U}={\EuScript U}(R, S)$
consisting of all $(0, 1)$-matrices of size $p$ by $q$ with row sum vector $R$ and column sum vector $S$.
If $R=(r,\ldots, r)$ and $S=(s,\ldots, s)$, for $r$ and $s$ positive integers, 
then ${\EuScript U}(R, S)$  is simply denoted by ${\EuScript U}(r, s)$. 
A  {\it sparse matrix} is a $(0,1)$-matrix in which most of the elements are zero.\\
If $A$ is an sparse matrix is important the $2$ by $2$ submatrices of $A$ of the types
$$ A_1= \begin{bmatrix} 1 & 0 \\ 0 & 1 \end{bmatrix} \quad {\rm and} \quad  
A_2= \begin{bmatrix} 0 & 1 \\ 1& 0\end{bmatrix}.$$
Now an {\it interchange} is a transformation of the elements of $A$ that changes a specified minor of type $A_1$ into type $A_2$ or else a minor of type $A_2$ into type $A_1$, and leaves all other elements of $A$ unaltered. Suppose that we apply to $A$ a finite number of interchanges. Then by the nature of the interchange operation, the resulting matrix $A^*$ has row sum vector $R$ and column sum vector $S$.
 For the convenience of the reader we state:
 
\begin{theorem}[Ryser]\label{Th-Ryser}
Let $A$ and $A^*$ be two $m$ by $n$ matrices composed of $0$'s and $1$'s,
possessing equal row sum vectors and equal column sum vectors. Then $A$ is
transformable into $A^*$ by a finite number of interchanges. 
\end{theorem} 

\begin{proof}
 See \cite[Theorem 3.1]{5}. 
\end{proof}

Following \cite{0.1}, let $X=\{ a_1,\ldots, a_n\}$ be a  set of $n$ elements. We call $X$ an $n$-{\it set}. Now $A_1,A_2,\ldots, A_m$ be $m$ not necessarily distinct subsets of the $n$-set $X$. We refer to this collection of subsets of an $n$-set  as a {\it  configuration of subsets}. 
Let $A=(a_{ij})$ a matrix, such that $a_{ij} = 1$ if $a_{ij} \in A_i$ and we set $a_{ij} = 0$ if $a_{ij} \not\in A_i$. The resulting $(0,1)$-matrix  $A=(a_{ij})$ with $i \in\{ 1,2,\ldots, m\}, j \in \{1,2,\ldots, n\}$ of size $m\times n$ is {\it the incidence matrix} for the configuration of subsets $A_1,A_2,\ldots, A_m$ of the $n$-set $X$.   Finally, an abstract system consisting of two types of objects and a single relationship between these types of objects is called {\it an incidence structure}.

We employ the following notation.
Let $m$ and $\ell$ be positive integers such that  $\ell< m$ as usual in the literature 
$C^m_{\ell}$ denotes binomial coefficient, we define {\it the index set} 
{ \small$$I(\ell,m)=\{\alpha=(\alpha_1,\ldots,\alpha_{\ell})\in {\mathbb N}^{\ell} :
1\leq\alpha_1<\cdots<\alpha_{\ell}\leq m\}.$$ }
 If $s\geq 1$ is a positive integer and  $\Sigma$ is a nonempty set,  we denote by $C_s(\Sigma)$ combinations of elements $\Sigma$ taken from $s$ ways. 
The elements of $C_s(\Sigma)$  are written as follows  $\big(P_{\alpha_1}, \ldots, P_{\alpha_s} \big)$ where $(\alpha_1, \ldots, \alpha_s)\in I(s, m)$ and we assume that $|\Sigma|=m$, the cardinality of $\Sigma$.  
Is important we notice that throughout this paper, when  $\alpha$  belongs a set of indexes $I(s,m)$, it is assumed that  belongs up a permutation that orders it properly.\\

We consider a symplectic  vector space and a subspace, a way to establish if this subspace is isotropic, 
is by means of a matrix and the Pl\"ucker coordinates, such matrix is given in \cite{2}.
In this paper for $n\geq 4$ arbitrary positive integer, $m=2k-2$, $2\leq k \leq r$ and $r=\lfloor\frac{n+2}{2} \rfloor$ partition is given in the whole $I(n-2, 2n)$ with sets $T_m$ called triangles and for each triangle $T_m$ defines a configuration and the incidence matrix of this configuration corresponds to a matrix  ${\EuScript L}_k$, we give a series of rules to form a matrix  $M={\EuScript L}_k\bigoplus\cdots\bigoplus {\EuScript L}_2$ which turns out to be a "canonical form" of the Lagrangian-Grassmannian variety $L(n, E)$ where $E$ is a simplectic vector space  
such that     $\dim E=2n$.\\

We define the matrix ${^t}A$  the {\it flip-transpose} of $A$ which flips $A$ across if skew-diagonal, 
so if $A=(a_{ij})_{mn}$, then 
${^t}A=(a_{n-j+1,m-i+1})_{nm}$.  A matrix we call {\it persymmetric} if ${^t}A = A$.
Let $A$ and $B$  be square matrices with coefficients in a field $\mathbb{F}$ of $n$ elements.  If there exist an invertible square matrix $P$ of order $n$ over $\mathbb{F}$  such that $B= P^{-1}AP$, then $A$ and $B$ are {\it similar}.

We denote by ${\mathbb J}$ the matrix, square, filled with ones and call it  {\it matrix of ones},  and $I^a_k$  denote the {\it anti-identity matrix} which an anti-diagonal square matrix,  filled with ones in the anti-diagonal. So we denote that 
$J_k={\mathbb J}-I_k$ where $I_k$ is the identity matrix of order  $k$ and note that the matrix $I^a_k$, with all anti-diagonal elements equal to 1, is the permutation matrix that reverses the order of elements of vectors.

A code is {\it formally self-dual} if the code and its dual have the same weight enumerator.  A code is 
{\it isodual} if it is equivalent to its dual. Clearly, any isodual code is formally self-dual.
 Binary codes are called even if provided have weights divisible by 2.\\
 Gleason's Theorem applies to even formally self-dual binary codes, states that
 if $C$ is such a code of length $n$, then 
 $$W_C=(x, y)=\sum_{i=0}^{\lfloor n/8\rfloor} a_ig_1(x, y)^{n/2-4i} g_2(x, y)^i$$ 
 where $g_1(x, y)=y^2+x^2$ and $g_2(x, y)=y^8+14x^4y^4+x^8$ are Gleason polynomials, see \cite{3.5}.

\section{Construction of sparse $(0,1)$-matrices} \label{Secc2}
We we use our algorithm to construct a structured, sparse, $(0,1)$-matrices $A_k^{\ell -2}$, in approximate lower diagonal form,  with $k$ ones in each row and $\ell$ ones in each column for $k$ and $\ell$ any positive integers.

Let $M$ be a matrix of order $n\times m$, the operation ${\EuScript O}(M)$ is to paste to the matrix $M$ at the bottom, the identity matrix $I_{m\times m}$, which generates a matrix $M (I_{m\times m})$  of order $(n+m)\times m$, that is ${\EuScript O}(M)=M(I_{m\times m})$. If we have a matrix vector $V=(M_1, \ldots, M_t)$  the operation 
${\EuScript O}(V)$ is the matrix vector $({\EuScript O}(M_1), \ldots, {\EuScript O}(M_t))$. 
Let $M_1$ and $M_2$ be two matrices of order $n_1\times m_1$ and   $n_2\times m_2$ respectively where $n_1\geq n_2$, 
{\it paste the matrix concatenatedly to the right}, side by side, is to get the matrix ${\EuScript P}(M_1, M_2)=M_1\bigsqcup M_2$  of order $n_1\times (m_1+m_2)$, where  $\sqcup$ means
joining together side-by-side and aligning the bottoms of the
corresponding identity matrices and filling the non-marked spaces on
the upper right blocks with zeroes. Finally denote by $\underline{k}$  the order matrix $1\times k$ filled with 1's.

\subsection{Algorithm to construct the matrix $A_k^{\ell-2}\in  {\EuScript U}(k, \ell)$}
In this subsection we present an algorithm to construct a structured sparse $(0,1)$-matrix 
in approximate lower triangular form,
with $k$ ones in each row and $\ell$ ones in each column 
where $k$ and $\ell$ any positive integers. 
 
\begin{center}
{\bf   Algorithm 1}\label{Algoritm1}
\end{center}
\hrulefill
\begin{description}
\item[{\rm Input: }] $k$ and $\ell$ arbitrary positive integers.
\item[{\rm Output: }] \; The matrix $A_k^{\ell-2}$ in approximate lower triangular form,
with $k$ ones in each row and  $\ell$ ones in each column.
\end{description}
\hrulefill
\begin{description}
\item[{\rm Step 1.}]\; Let $V=(\underline{k}, \underline{k-1}, \ldots, \underline{2}, \underline{1})$ be matrix vector.
\item[{\rm Step 2.}] \;We apply the operation ${\EuScript O}(V)$ to the matrix vector given in Step 1 and  we obtain   
the matrix vector $({\EuScript O}(\underline{k}), {\EuScript O}(\underline{k-1}), \ldots, {\EuScript O}(\underline{2}), {\EuScript O}(\underline{1}))$.
\item[{\rm Step 3.}]\; Put ${\EuScript P}({\EuScript O}(\underline{k}), {\EuScript O}(\underline{k-1}), \ldots, {\EuScript O}(\underline{2}), {\EuScript O}(\underline{1}))$ we generate a matrix that we denote by $A_k^0$.
\item[{\rm Step 4.}]\; Now consider the matrix vector \\ $V=(A_k^0, A_{k-1}^0,\ldots, A_2^0,A_1^0)$.
\item[{\rm Step 5.}]\; Returne to the Steps 2,  3,  4  with matrix vector $V=(A_k^0, A_{k-1}^0,\ldots, A_2^0,A_1^0)$
to build a matrix $A_k^1$  given by ${\EuScript P}({\EuScript O}(A_k^0), {\EuScript O}(A_{k-1}^0),\ldots, {\EuScript O}(A_2^0), {\EuScript O}(A_1^0))$.
\item[{\rm Step 6.}] \;The algorithm ends when $k=\ell-2$.
\end{description}
\hrulefill

\begin{definition}\label{Ass-Matrix}
 Let  $k$ and $\ell$ be arbitrary positive integers,  we define the matrices  $A_k^{\ell-2}$  inductively as
 $A_k^0:={\EuScript P}({\EuScript O}(\underline{k}), {\EuScript O}(\underline{k-1}), \ldots, {\EuScript O}(\underline{2}), {\EuScript O}(\underline{1}))$ and \\
 $A_k^1:= {\EuScript P}({\EuScript O}(A_k^0), {\EuScript O}(A_{k-1}^0),\ldots, {\EuScript O}(A_2^0), {\EuScript O}(A_1^0))$, that is 
  $$ A_k^1=A_k^0 \left(I_{C^{k+1}_2}
    \right)\sqcup A_{k-1}^0 \left( I_{C^{k}_2} \right)\sqcup \cdots \sqcup
    A_2^0\left( I_{C^3_2} \right)\sqcup A_1^0 \left( I_{C^2_2} \right).$$
 With the previous notation we define the following matrices
 \begin{align*}
    A_k^2&=A_k^1\left( I_{C^{k+2}_3}
    \right) \sqcup A_{k-1}^1 \left( I_{C^{k+1}_3} \right) \sqcup \cdots \sqcup
    A_2^1\left( I_{C^4_3} \right) \sqcup A_1^1\left( I_{C^3_3} \right),\\
    A_k^3&=A_k^2\left( I_{C^{k+3}_4}
    \right) \sqcup A_{k-1}^2 \left( I_{C^{k+2}_4} \right)\sqcup \cdots \sqcup
    A_2^2\left( I_{C^5_4} \right)\sqcup A_1^2 \left( I_{C^4_4}\right),\\
    \vdots & \quad\quad\quad\quad\quad\quad\quad \quad\quad\quad\vdots \\
  A_k^{\ell-2} & =A_k^{\ell-3}\left( I_{C^{k+\ell-2}_{\ell-1}}
\right)\sqcup  \cdots \sqcup A_2^{\ell-3}\left(
I_{C^{\ell}_{\ell-1}} \right) \sqcup A_1^{\ell-3}\left(
I_{C^{\ell-1}_{\ell-1}} \right). 
  \end{align*}
 \end{definition}

\subsection{Properties of $A_k^{\ell -2}$}\label{Prop}

In this subsection we present some properties of the family of matrices $A_k^{\ell -2}$.

\begin{theorem}\label{exis}
Let $k$ and $\ell$ be any positive integers. Then $A_k^{\ell -2}$ is a sparse $(0, 1)$-matrix
 with $k$ ones in each row and $\ell$ ones in each column, and is of order  
$C_{\ell-1}^{k+\ell-1}\times C_{\ell}^{k+\ell-1}$. 
\end{theorem}

\begin{proof}
We will show that  the matrix $A_k^{\ell-2}$ given above satisfies the theorem conditions.
The proof is by induction on $\ell$, number of ones in each column and k  a fixed arbitrary positive integer.
If $\ell=2$ evidently $A_k^0\in {\EuScript U}(k, 2)$ for  $k$ any positive integer.
Now by induction hypothesis suppose that $A_k^{\ell-3}\in {\EuScript U}(k, \ell-1)$.  
It is easy to see that the matrix $A_k^{\ell-3}\left( I_{C^{k+\ell-2}_{\ell-1}}
\right)$ has $\ell$ ones in each column and $k$  ones in a part of the matrix rows and exactly 1 in the rest of the rows. \\
Then the matrix $A_k^{\ell-3}\left( I_{C^{k+\ell-2}_{\ell-1}}
\right)\sqcup A_{k-1}^{\ell-3}\left( I_{C^{k+\ell-3}_{\ell-1}}
\right)$ has $\ell$ ones in each column and a greater number of rows with $k$  ones and the rest of rows with exactly two ones.
 Continuing in this way we obtain  the matrix  $$A_k^{\ell-3}\left( I_{C^{k+\ell-2}_{\ell-1}}
\right)\sqcup A_{k-1}^{\ell-3}\left( I_{C^{k+\ell-3}_{\ell-1}}
\right)   \sqcup  \cdots \sqcup A_2^{\ell-3}\left(
I_{C^{\ell}_{\ell-1}} \right)$$ 
this matrix has $\ell$ ones in each column and exactly $\ell$ rows with exactly $k-1$ ones in each row. Then  the matrix $$A_k^{\ell-2}  =A_k^{\ell-3}\left( I_{C^{k+\ell-2}_{\ell-1}}
\right)\sqcup  \cdots \sqcup A_2^{\ell-3}\left(
I_{C^{\ell}_{\ell-1}} \right) \sqcup A_1^{\ell-3}\left(
I_{C^{\ell-1}_{\ell-1}} \right) $$ has $k$ ones in each row and $\ell$ ones in  each column. 
 To find the order of the matrix, we observe that its number of rows is
$$1+C^k_1+C^{k+1}_{2}+\cdots+C^{k+\ell-3}_{\ell-2}+C^{k+\ell-2}_{\ell-1}=C^{k+\ell-1}_{\ell-1},$$
and its number of columns is
$$C^{k+\ell-2}_{\ell-1}+C^{k+\ell-3}_{\ell-1}+\cdots+C^{\ell}_{\ell-1}+C^{\ell-1}_{\ell-1}=C^{k+\ell-1}_{\ell}.$$
Moreover  $A_k^{\ell-2}$  is a sparse $(0,1)$-matrix, indeed the quotient formed by number of ones of the matrix
between the order is $\frac{kC^{k+\ell-1}_{\ell-1}}{C^{k+\ell-1}_{\ell-1}\times C^{k+\ell-1}_{\ell}}$ and tends to zero when $k$ and $\ell$ tends to infinite.
\end{proof}

\begin{corollary}
The matrix $A_k^{\ell-2}$ is square if and only if  $\ell=k$.
\end{corollary}

\begin{proof}
By Theorem  \ref{exis} $A_k^{\ell-2}$ has orden $C^{k+\ell-1}_{\ell-1}\times C^{k+\ell-1}_{\ell}$. 
Hence  $A_k^{\ell-2}$  is square if and only if  $C^{k+\ell-1}_{\ell-1}= C^{k+\ell-1}_{\ell}$ if and only if  $k+\ell-1=\ell-1+\ell$ if and only if  $k=\ell$.
\end{proof}

\begin{theorem}\label{trans}
Let $k$ and  $\ell$ be  arbitrary positive integers. Then  $\big(A_k^{\ell-2}\big)^t$ is transformable in $A_{\ell}^{k-2}$.
\end{theorem}

\begin{proof}
We notice that the matrix $A_k^{\ell-2}\in {\EuScript U}(k, \ell)$  has orden $C_{\ell-1}^{k+\ell-1}\times C_{\ell}^{k+\ell-1}$ and the transposed matrix  
$\big(A_k^{\ell-2}\big)^t\in {\EuScript U}(\ell, k)$ has orden $C_{\ell}^{\ell+ k-1}\times C_{\ell -1}^{\ell + k-1}$. On the other hand 
 $A_{\ell}^{k-2}\in {\EuScript U}(\ell, k)$  has orden $C_{k-1}^{\ell+ k-1}\times C_{k}^{\ell + k-1}$ and also 
 $C_{\ell}^{\ell+ k-1}\times C_{\ell -1}^{\ell + k-1}=C_{k-1}^{\ell+ k-1}\times C_{k}^{\ell + k-1}$. Therefore  $\big(A_k^{\ell-2}\big)^t$ and  $A_{\ell}^{k-2}$ have
the same order and $\big(A_k^{\ell-2}\big)^t$, $ A_{\ell}^{k-2} \in {\EuScript U}(\ell, k)$, then by Theorem  \ref{Th-Ryser},  $\big(A_k^{\ell-2}\big)^t$ is transformable in 
$A_{\ell}^{k-2} $ with a finite number of interchanges. 
\end{proof}

\begin{definition}\label{Def-Mat-Fragm}
Let $n$, $m$, $r$ and $s$ positive integers such that $m> n$ and $n-s=m-r>0$. Let $A$ a matrix of order $n \times m$, $B$ a matrix of order $s\times (m-r)$, $C$
 a square matrix of order $ (m-r)\times (n-s)$, 
 and  $D$ a matrix of order $(n-s)\times r$, we will call $A$ a  {\it fragmented matrix} if
 $A= B(C)\sqcup D$ where  $B(C)$  is to paste to the matrix $B$ at bottom, the matrix $C$ and $\sqcup$ 
 means joining together side-by-side and aligning the bottoms  the  matrix $D$ and filling the non-marked spaces on the upper right blocks with 0's.
 \end{definition}

\begin{theorem}\label{Fragment-11}
Let $k$ and $\ell$ be any positive integers. Then the matrix $A_k^{\ell-2}$ is a fragmented matrix as follows 
$A_k^{\ell-2} = A_k^{\ell-3}(I_{C_{\ell -1}^{k+\ell -2}}) \sqcup A_{k-1}^{\ell-2}.$ 
\end{theorem}

\begin{proof}
By Algorithm \ref{Algoritm1} and Theorem \ref{exis} we have that $A_k^{\ell-2}$ is constructed with a recursive algorithm.
That is
$$ A_k^{\ell-2}  =A_k^{\ell-3}\left( I_{C^{k+\ell-2}_{\ell-1}}
\right)\sqcup  A_{k-1}^{\ell-3}\left(
I_{C^{ k+\ell-3}_{\ell-1}} \right) \sqcup  \cdots \sqcup A_1^{\ell-3}\left(
I_{C^{\ell-1}_{\ell-1}} \right).$$
By definition we have
$$A_{k-1}^{\ell-2}=A_{k-1}^{\ell-3}\left(
I_{C^{ k+\ell-3}_{\ell-1}} \right) \sqcup  \cdots \sqcup A_1^{\ell-3}\left(
I_{C^{\ell-1}_{\ell-1}} \right). $$
Then
$$A_k^{\ell-2} = A_k^{\ell-3}(I_{C_{\ell -1}^{k+\ell -2}}) \sqcup A_{k-1}^{\ell-2}.$$
With this we show the fragmented of $A_k^{\ell-2}$. 
\end{proof}

\begin{corollary}\label{app-lower-diag}
Let $k$ and $\ell$ be any positive integers. Then $A_k^{\ell -2}$ is a  $(0, 1)$-matrix
approximate lower diagonal form.
\end{corollary}
\begin{proof}
From Theorem \ref{Fragment-11} we have that the matrix   $A_k^{\ell-2} = A_k^{\ell-3}(I_{C_{\ell -1}^{k+\ell -2}}) \sqcup A_{k-1}^{\ell-2}$ is a fragmented matrix, so
$I_{C_{\ell -1}^{k+\ell -2}}$ is a submatrix, in the lower left corner,  and by definition
 $A_k^{\ell -2}$ is a  matrix approximate lower diagonal form.
\end{proof}

\begin{theorem}\label{propT1} Let $k$ and $\ell$ any positive integers. Then 
$${^t}A_{k}^{\ell -2} = A_{\ell}^{k-2}.$$
\end{theorem}

\begin{proof}
First it is direct to see that both matrices have the same order.  We show the statement using induction on $k$ and $\ell$.
For $k$ arbitrary positive integer and $\ell =1$, we have 
$$ ^t(A_k^{1-2})= \begin{bmatrix} 1\\ 1 \\ \vdots \\ 1 \end{bmatrix} = A_1^{k-2}.$$
If $\ell$ arbitrary positive integer and $k = 1$, we have
$$ ^t(A_1^{\ell-2})= ( 1, 1,  \ldots, 1) = A_\ell^{1-2}.$$
Assume the true result for  all $k' < k$ and $\ell' < \ell$, i.e., $ {^t}A_{k'}^{\ell -2} = A_{\ell}^{k'-2}$  and ${^t}A_{k}^{\ell' -2} = A_{\ell'}^{k-2}$. Now for $k$ and $\ell$ positive integers we have
\begin{align*}
A_k^{\ell-2} &=  A_k^{\ell-3}(I_{C_{\ell -1}^{k+\ell -2}}) \sqcup A_{k-1}^{\ell-2} \\
 {^t}(A_k^{\ell-2}) & =  {^t}\left( A_k^{\ell-3}(I_{C_{\ell -1}^{k+\ell -2}}) \sqcup A_{k-1}^{\ell-2}\right)\\
 {^t}(A_k^{\ell-2}) & =  {^t}(A_{k-1}^{\ell-2})(I_{C_{\ell -1}^{k+\ell -2}}) \sqcup {^t}(A_{k}^{\ell-3})
 \end{align*}
 by induction we obtain
 $$
 {^t}(A_k^{\ell-2})  =  (A_\ell^{k-3})(I_{C_{\ell -1}^{k+\ell -2}}) \sqcup A_{\ell-1}^{k-2} =  A_\ell^{k-2}.
 $$ 
\end{proof}

\begin{corollary}\label{cuaPe}
If $\ell = k$, then $A_k^{k-2}$ is persymmetric and the matrix $A_k^{k-2}$ is similar to $(A_k^{k-2})^t$, up permutation of rows and columns.
\end{corollary}
\begin{proof}
From Theorem  \ref{propT1} we have ${^t}A_{k}^{k -2} = A_{k}^{k-2}$ and so we prove that $A_k^{k-2}$ is persymmetric. Now it is easy to see $I^a_{C_{k-1}^{2k-2}} (^tA_k^{k-2})I^a_{C_{k-1}^{2k-2}} = (A_k^{k-2})^t $ where $I^a_{C_{k-1}^{2k-2}}$ denote the anti-identity matrix. With this we show the similarity of $A_k^{k-2}=^t(A_k^{k-2})$ and $(A_k^{k-2})^t$.
\end{proof}

\section{Construction of $B_k^{\ell-2}$ dense $(0,1)$-matrices} \label{Secc3}

Let  $ {\mathbb J}$ matrix full of ones, and let $M$ be a matrix of order $n\times m$, the operation ${\EuScript O}(M)$ is to paste to the matrix $M$ at the bottom, the matrix $J_{m\times m}= {\mathbb J}-I_{m\times m}$  which generates a matrix $M (J_{m\times m})$  of order $(n+m)\times m$, that is ${\EuScript O}(M)=M(J_{m\times m})$. If we have a matrix vector $V=(M_1, \ldots, M_t)$  the operation 
${\EuScript O}(V)$ is the matrix vector $({\EuScript O}(M_1), \ldots, {\EuScript O}(M_t))$. 
Let $M_1$ and $M_2$ be two matrices of order $n_1\times m_1$ and   $n_2\times m_2$ respectively 
where $n_1\geq n_2$, 
{\it paste the matrix concatenatedly to the right}, side by side, is to get the matrix ${\EuScript P}(M_1, M_2)=M_1\bigsqcup M_2$  of order $n_1\times (m_1+m_2)$, where  $\sqcup$ means
joining together side-by-side and aligning the bottoms of the
corresponding matrices and filling the non-marked spaces on
the upper right blocks with $1$'s. Finally denote by $\overline{k}=(0, \ldots, 0)$  the order matrix $1\times k$.
 
\subsection{Algorithm to construct the matrix $B_k^{\ell-2}\in  {\EuScript V}(C_{\ell-1}^{k+\ell-1}-k, C_{\ell-1}^{k+\ell-1}-\ell)$}
In this part we adapt the algorithm \ref{Algoritm1} to construct a structured dense $(0,1)$-matrix 
such that $A_k^{\ell-2}+B_k^{\ell-2}={\mathbb J}$. 
\begin{center}
{\bf   Algorithm 2}\label{Algoritm 2} 
\end{center}
\hrulefill
\begin{description}
\item[{\rm Input: }] $k$ and $\ell$ arbitrary positive integers.
\item[{\rm Output: }] \; The dense matrix $B_k^{\ell-2}$ such that $A_k^{\ell-2}+B_k^{\ell-2}={\mathbb J}$.
\end{description}
\hrulefill
\begin{description}
\item[{\rm Step 1.}]\; Let $V=(\overline{k}, \overline{k-1}, \ldots, \overline{2}, \overline{1})$ be matrix vector.
\item[{\rm Step 2.}] \;We apply the operation ${\EuScript O}(V)$ to the matrix vector given in Step 1 and  we obtain   
the matrix vector $({\EuScript O}(\overline{k}), {\EuScript O}(\overline{k-1}), \ldots, {\EuScript O}(\overline{2}), {\EuScript O}(\overline{1}))$.
\item[{\rm Step 3.}]\; Put ${\EuScript P}({\EuScript O}(\overline{k}), {\EuScript O}(\overline{k-1}), \ldots, {\EuScript O}(\overline{2}), {\EuScript O}(\overline{1}))$ we generate a matrix that we denote by $B_k^0$.
\item[{\rm Step 4.}]\; Now consider the matrix vector \\ $V=(B_k^0, B_{k-1}^0,\ldots, B_2^0, B_1^0)$.
\item[{\rm Step 5.}]\; Returne to the Steps 2,  3,  4  with matrix vector $V=(B_k^0, B_{k-1}^0,\ldots, B_2^0, B_1^0)$
to build a matrix $B_k^1$  given by ${\EuScript P}({\EuScript O}(B_k^0), {\EuScript O}(B_{k-1}^0),\ldots, {\EuScript O}(B_2^0), {\EuScript O}(B_1^0))$.
\item[{\rm Step 6.}] \;The algorithm ends when $k=\ell-2$.
\end{description}
\hrulefill

\begin{definition}\label{Ass22-Matrix}
 Let  $k$ and $\ell$ be arbitrary positive integers,  we define the matrices  $B_k^{\ell-2}$  inductively as
 $B_k^0:={\EuScript P}({\EuScript O}(\underline{k}), {\EuScript O}(\underline{k-1}), \ldots, {\EuScript O}(\underline{2}), {\EuScript O}(\underline{1}))$ and \\
 $B_k^1:= {\EuScript P}({\EuScript O}(B_k^0), {\EuScript O}(B_{k-1}^0),\ldots, {\EuScript O}(B_2^0), {\EuScript O}(B_1^0))$, that is 
  $$ B_k^1=B_k^0 \left(J_{C^{k+1}_2}
    \right)\sqcup B_{k-1}^0 \left( J_{C^{k}_2} \right)\sqcup \cdots \sqcup
    B_2^0\left( J_{C^3_2} \right)\sqcup B_1^0 \left( J_{C^2_2} \right).$$
 With the previous notation we define the following matrices
 \begin{align*}
    B_k^2&=B_k^1\left( J_{C^{k+2}_3}
    \right) \sqcup B_{k-1}^1 \left( J_{C^{k+1}_3} \right) \sqcup \cdots \sqcup
    B_2^1\left( J_{C^4_3} \right) \sqcup B_1^1\left( J_{C^3_3} \right),\\
    B_k^3&=B_k^2\left( J_{C^{k+3}_4}
    \right) \sqcup B_{k-1}^2 \left( J_{C^{k+2}_4} \right)\sqcup \cdots \sqcup
    B_2^2\left( J_{C^5_4} \right)\sqcup B_1^2 \left( J_{C^4_4}\right),\\
    \vdots & \quad\quad\quad\quad\quad\quad\quad \quad\quad\quad\vdots \\
  B_k^{\ell-2} & =B_k^{\ell-3}\left( J_{C^{k+\ell-2}_{\ell-1}}
\right)\sqcup  \cdots \sqcup B_2^{\ell-3}\left(
J_{C^{\ell}_{\ell-1}} \right) \sqcup B_1^{\ell-3}\left(
J_{C^{\ell-1}_{\ell-1}} \right). 
  \end{align*}
 \end{definition}

\subsection{Properties of $B_k^{\ell -2}$}\label{Prop}

In this subsection we present some properties of the family of matrices $B_k^{\ell -2}$.

\begin{theorem}\label{exis2}
The $(0, 1)$-matrix $B_k^{\ell-2}$ is dense, of the same order as $A_k^{\ell-2}$  such that $B_k^{\ell-2}={\mathbb J}-A_k^{\ell-2}$.
\end{theorem}

\begin{proof}
Is enough to prove  $B_k^{\ell-2}={\mathbb J}-A_k^{\ell-2}$. We use induction on $\ell$ and $k$.
If $\ell=2$ and $k$ any arbitrary positive integer it is easy to see that  $B_k^{2-2}={\mathbb J}-A_k^{2-2}$. 
Also if $k=2$ and $\ell$ any arbitrary positive integer,   then  it is easy to see that  $B_2^{\ell-2}={\mathbb J}-A_2^{\ell-2}$.
Now the induction hypothesis is that for each  $k^{'}<k$ and $\ell^{'}< \ell$ arbitrary positive integer, 
then    $B_{k^{'}}^{{\ell^{'}}-2}={\mathbb J}-A_{k^{'}}^{\ell^{'}-2}$ where $1\leq  \ell^{'} \leq \ell$.  
 We have that  the matrix 
 \begin{align*}
 {\mathbb J}-A_k^{\ell-2} &= {\mathbb J}- \big[A_k^{\ell-3}\left( I_{C^{k+\ell-2}_{\ell-1}}
\right)\sqcup  \cdots \sqcup  A_1^{\ell-3}\left(
I_{C^{\ell-1}_{\ell-1}} \right)\big]\\
 & =[{\mathbb J}-A_k^{\ell-3}\left( I_{C^{k+\ell-2}_{\ell-1}}
\right)] \sqcup  \cdots \sqcup [ {\mathbb J}-A_1^{\ell-3}\left(
I_{C^{\ell-1}_{\ell-1}} \right)] \\
&=B_k^{\ell-3}\left( I_{C^{k+\ell-2}_{\ell-1}}
\right)\sqcup  \cdots \sqcup  B_1^{\ell-3}\left(
I_{C^{\ell-1}_{\ell-1}} \right) \\
&=B_k^{\ell-2}.
\end{align*}
\end{proof}

\begin{corollary}\label{Pers-dense}
The matrix $B_k^{\ell-2}$ is square if and only if  $\ell=k$.
\end{corollary}

\begin{proof}
By Theorem  \ref{exis2} $B_k^{\ell-2}$ has orden $C^{k+\ell-1}_{\ell-1}\times C^{k+\ell-1}_{\ell}$. 
Hence  $B_k^{\ell-2}$  is square if and only if  $C^{k+\ell-1}_{\ell-1}= C^{k+\ell-1}_{\ell}$ if and only if  $k+\ell-1=\ell-1+\ell$ if and only if  $k=\ell$.
\end{proof}

\begin{theorem}\label{denseFragmente}
Let $k$ and $\ell$ be any positive integers. Then the matrix $B_k^{\ell-2}$  is a fragmented matrix. 
\end{theorem}

\begin{proof}
By Theorem \ref{Fragment-11} we have 
\begin{align*}
B_k^{\ell-2}&={\mathbb J}-A_k^{\ell-2}\\
&={\mathbb J}-\big[A_k^{\ell-3}(I_{C_{\ell-1}^{k+\ell-2}})\sqcup A_{k-1}^{\ell-2}\big]\\
&=\big[{\mathbb J}-A_k^{\ell-3}(I_{C_{\ell-1}^{k+\ell-2}})\big] \sqcup [{\mathbb J}-A_{k-1}^{\ell-2}]\\
&=B_k^{\ell-3}(J_{C_{\ell-1}^{k+\ell-2}})\sqcup B_{k-1}^{\ell-2}.
\end{align*}
\end{proof}

\begin{theorem}\label{propT} 
Let $k$ and $\ell$ any positive integers. Then 
$${^t}B_{k}^{\ell -2} = B_{\ell}^{k-2}.$$
\end{theorem}

\begin{proof}
In effect $B_k^{\ell-2}={\mathbb J}-A_k^{\ell-2}$, see Theorem \ref{exis2}, then
$$
^tB_k^{\ell-2} =^t\big({\mathbb J}-A_k^{\ell-2}\big)
=^t{\mathbb J}-^tA_k^{\ell-2}
={\mathbb J}-A_{\ell}^{k-2}
=B_{\ell}^{k-2}.$$
\end{proof}

\begin{corollary}\label{cuaPe2}
If $\ell = k$, then $B_k^{k-2}$ is persymmetric and the matrix $B_k^{k-2}$ is similar to $(B_k^{k-2})^t$, by permutation of rows and columns. 
\end{corollary}

\begin{proof}
By Corollary \ref{Pers-dense}  it easily follows that $B_k^{k-2}$ is persymmetric.  
We can easily see $I^a_{C_{k-1}^{2k-2}} (B_k^{k-2})I^a_{C_{k-1}^{2k-2}} = I^a_{C_{k-1}^{2k-2}} (^tB_k^{k-2})I^a_{C_{k-1}^{2k-2}} =(B_k^{k-2})^t $ where $I^a_{C_{k-1}^{2k-2}}$
denote the anti-identity matrix. 
Which proves the theorem.
\end{proof}
\section{ $ A_k^{k-3}$  incidence matrix of the configuration}\label{Secc4}
In this section, we show that the matrix $A_k^{k-3} \in {\EuScript U}(k, k-1)$, is
incidence matrix of the configuration of triangle-sets  $T_m$ of the set of indexes.  Furthermore, we proof that a sparse $(0,1)$-matrix $M$  it is expressed as a direct sum of matrices $A_k^{k-3}$ see Theorems \ref{thmPeven} and \ref{thmPodd}.
 
Throughout this section we will denote by $m$ and $\ell$ positive integers such that  $\ell< m$ and consider the index set   
$$I(\ell, m)=\{\alpha=(\alpha_1,\ldots, \alpha_{\ell})\in {\mathbb N}^{\ell}: 1\leq\alpha_1<\cdots <\alpha_{\ell}\leq m \}.$$ 
Now let $m\in{\mathbb N}$ be an even integer and define: 
 \begin{align*}
P_i&=(i,2m-i+1)\in {\mathbb N}^2,\quad\text{for $1\leq i\leq m$,}\\
\Sigma_s&=\{P_1,\ldots,P_s\}\subset {\mathbb N}^2, \quad\text{for $1\leq s\leq m$,}\\
\Sigma_{[s,m]}&=\{P_{s+1},\ldots,P_m\}\subset {\mathbb N}^2.
\end{align*}
As in preliminaries, we denote by $C_{(m-6)/2}(\Sigma_m)$  the 
combinations of elements $\Sigma_m $ taken from $(m-6)/2$ forms.
Now for all  $(\alpha(1),\ldots,\alpha((m-6)/2))\in I((m-6)/2,m-2)$ with $m\geq8$ even integer,
we define the set  
$$T_{(\alpha(1),\alpha(2),\ldots,\alpha(\frac{m-8}{2}),\alpha(\frac{m-6}{2}))} := T_{\alpha(1,2,\ldots, \frac{m-6}{2})},$$
 as the cartesian product 
$$T_{\alpha(1,2,\ldots, \frac{m-6}{2})}=(P_{\alpha(1)},\ldots,P_{\alpha((m-8)/2)},
P_{\alpha((m-6)/2)}) \times C_2(\Sigma_{[\alpha((m-6)/2),m]}).$$

We are interested in putting all the necessary notation to establish our main result of this subsection. Following this end
{\footnotesize
$${\begin{array}{l}
    T_{1}^1 =T_{\alpha(1,2,\ldots,\frac{m-8}{2},\frac{m-6}{2})}\\
    T_{2}^1 =\bigcup_{i=0}^{\frac{m+2}{2}}T_{\alpha(1, 2,\ldots,
        \frac{m-8}{2}, \frac{m-6}{2}+i)}\\
    T_{3}^1 =T_{2}^1\cup\left(\bigcup_{i=0}^{\frac{m}{2}}T_{\alpha(1, 3,\ldots,
        \frac{m-6}{2}, \frac{m-4}{2}+i)}\right)\\
    T_{4}^{1} =T_{3}^1\cup \left(\bigcup_{i=0}^{\frac{m-2}{2}}T_{\alpha(1,
        4,\ldots,\frac{m-4}{2},\frac{m-2}{2}+i)}\right)\\
    \;\;\,\quad \vdots  \\
    T_{\frac{m+4}{2}}^{1} =T_{\frac{m+2}{2}}^1\cup \left(\bigcup_{i=0}^{1}T_{\alpha(1, \frac{m+4}{2},\ldots, m-3+i)}\right)\\
    T_{\frac{m+6}{2}}^{1} =T_{\frac{m+4}{2}}^1\cup T_{\alpha(1,
        \frac{m+6}{2},\ldots,m-3,m-2)}
    \end{array}}$$ 
   
$${\begin{array}{l}
    T_{2}^2 =T_{\alpha(2,3,\ldots,\frac{m-6}{2},\frac{m-4}{2})}\\
    T_{3}^2 =\bigcup_{i=0}^{\frac{m}{2}}T_{\alpha(2,3,\ldots, \frac{m-6}{2},
        \frac{m-4}{2}+i)}\\
    T_{4}^2 =T_{3}^2\cup\left(\bigcup_{i=0}^{\frac{m-2}{2}}T_{\alpha(2,4,\ldots,
        \frac{m-4}{2}, \frac{m-2}{2}+i)}\right)\\
    T_{5}^{2} =T_{4}^2\cup
    \left(\bigcup_{i=0}^{\frac{m-4}{2}}T_{\alpha(2,5,\ldots,\frac{m-2}{2},
        \frac{m}{2}+i)}\right)\\
    \;\;\,\quad \vdots \\
    T_{\frac{m+4}{2}}^{2} =T_{\frac{m+2}{2}}^2\cup\left(\bigcup_{i=0}^{1}T_{\alpha(2, \frac{m+4}{2},\ldots,m-3+i)}\right)\\
    T_{\frac{m+6}{2}}^{2} =T_{\frac{m+4}{2}}^2\cup T_{\alpha(2,
        \frac{m+6}{2},\ldots,m-2)}
    \end{array}}
 $$
 %
$${\begin{array}{l}
    T_{3}^3 = T_{\alpha(3,4,\ldots,\frac{m-4}{2},\frac{m-2}{2})}\\
    T_{4}^3 = \bigcup_{i=0}^{\frac{m-2}{2}}T_{\alpha(3,4,\ldots, \frac{m-4}{2},\frac{m-2}{2}+i)}\\
    T_{5}^3 = T_{4}^3\cup\left(\bigcup_{i=0}^{\frac{m-4}{2}}T_{\alpha(3,5,\ldots,
        \frac{m-2}{2},
        \frac{m}{2}+i)}\right)\\
    T_{6}^{3} = T_{5}^3\cup
    \left(\bigcup_{i=0}^{\frac{m-6}{2}}T_{\alpha(3,6,\ldots,\frac{m}{2},
        \frac{m+2}{2}+i)} \right)\\
    \;\;\,\quad \vdots    \\
    T_{\frac{m+4}{2}}^{3} = T_{\frac{m+2}{2}}^3\cup\left(\bigcup_{i=0}^{1}T_{\alpha(3, \frac{m+4}{2}, \ldots,m-3+i)}
    \right)\\
    T_{\frac{m+6}{2}}^{3} = T_{\frac{m+4}{2}}^3\cup T_{\alpha(3,
        \frac{m+6}{2},\ldots,m-2)}
    \end{array}}$$ 
    $$ \vdots$$
$${\begin{array}{l}
    T_{\frac{m+2}{2}}^{\frac{m+2}{2}}=T_{\alpha(\frac{m+2}{2},\frac{m+4}{2},\ldots,m-4,m-2)}\\
    T_{\frac{m+4}{2}}^{\frac{m+2}{2}}=\bigcup_{i=0}^1T_{\alpha(\frac{m+2}{2},\frac{m+4}{2},\ldots,m-4,m-3+i)}\\
    T_{\frac{m+6}{2}}^{\frac{m+2}{2}}=T_{\frac{m+4}{2}}^{\frac{m+2}{2}}\cup T_{\alpha(\frac{m+2}{2},\frac{m+6}{2},\ldots,m-4,m-2)}\\
    T_{\frac{m+6}{2}}^{\frac{m+4}{2}}=T_{\alpha(\frac{m+4}{2},\frac{m+6}{2},\ldots,m-3,m-2)}.
    \end{array}}$$}
For $m\geq 8$  we define the set 
\begin{equation}\label{eqTriang}
T_{m}=T_{\frac{m+6}{2}}^{1}\cup T_{\frac{m+6}{2}}^{2}\cup
\cdots \cup  T_{\frac{m+6}{2}}^{\frac{m+2}{2}} \cup
T_{\frac{m+6}{2}}^{\frac{m+4}{2}}.
\end{equation}

\begin{remark}
It is easy to see that  $T_m\subseteq C_{\frac{m-2}{2}}(\Sigma_m)$ where $m\geq 8$.
\end{remark}

\begin{definition}\label{TrianG}
If $m$ is even integer we define the set 
$$
\overline{T}_m=\begin{cases}
C_0(\Sigma_2) & \text{if $m=2$}, \\
C_1(\Sigma_4) & \text{if $m=4$}, \\
C_2(\Sigma_6) & \text{if $m=6$}, \\ T_m
& \text{if $m\geq 8$}, \\
\end{cases}$$
we call $\overline {T}_m$ a { \it triangulation} of the set $C_{\frac{m-2}{2}}(\Sigma_m)$.
 If $C_{\frac{m-2}{2}}(\Sigma_m)=\overline{T}_m$ then we say that $C_{\frac{m-2}{2}}(\Sigma_m)$
 is a { \it triangulated set}.
\end{definition}

\subsection{ An incidence function  $\varphi^m$}\label{phifuntion}
For $J$ a set not empty we denote by $$\{ 0, 1\}^J= \{ (a_j)_{j\in J} : a_j=0 \;{\rm or }\; a_j=1\}.$$
Let  $m$ be even integer and let $\varphi^m:T_m\rightarrow \{0,1\}^{I(m/2,m)}$  the
function given by
 \begin{equation}\label{Sec-Mm}
(P_{\alpha(1)},\ldots,P_{\alpha((m-2)/2)})\mapsto \big( \varphi_{\beta}^m(P_{\alpha(1)},\ldots,P_{\alpha((m-2)/2)}) \big)_{\beta\in I(m/2,m)\ },
\end{equation}    
where $\varphi_{\beta}^m(P_{\alpha(1)},\ldots,P_{\alpha((m-2)/2)}) \big)=$ 
 $$=\begin{cases}
1  & \text{if $\beta\in\{(P_{\alpha(1)},\ldots,P_{\alpha(\frac{m-2}{2})},P_i) : |\text{supp}(P_{\alpha(1)},\ldots,P_{\alpha(\frac{m-2}{2})},P_i )| = m\} $}, \\
0& \text{otherwise}.
\end{cases}$$
with 
$\text{supp}\{(P_{\alpha(1)},\ldots,P_{\alpha((m-2)/2)},P_i) \}
=\{ \alpha(1), 2n-\alpha(1)-1, \ldots, \alpha((m-2)/2), 2n-\alpha((m-2)/2)-1, \alpha(i), 2n-\alpha(i)-1 \}$.\\
Clearly
$$ |\{(P_{\alpha(1)},\ldots,P_{\alpha((m-2)/2)},P_i) : |{\rm supp}
(P_{\alpha(1)},  \ldots, P_{\alpha((m-2)/2)}, P_i| = m \}| $$ $$= \frac{m+2}{2},$$  and we define $k= \frac{m+2}{2}.$ Then the weight of the vector
$\varphi^m((P_{\alpha(1)},\ldots, P_{\alpha(\frac{m-2}{2})}))$ is $k$ for all  $(P_{\alpha(1)},\ldots, P_{\alpha(\frac{m-2}{2})}) \in T_m.$

\begin{lemma}\label{Funt-Inj} 
Let  $m\geq 2$ be even integer. Then the function $\varphi^m$  is injective.
\end{lemma}

\begin{proof}
Let  $(P_{\alpha(1)},\ldots,P_{\alpha((m-2)/2)}),  (P^{'}_{\alpha(1)},\ldots,P^{'}_{\alpha((m-2)/2)})\in T_m$ such that 

 $\varphi^m((P_{\alpha(1)},\ldots,P_{\alpha((m-2)/2)}))=\varphi^m((P^{'}_{\alpha(1)},\ldots,P^{'}_{\alpha((m-2)/2)}))$.
Then \\
$ \varphi_{\beta}^m(P_{\alpha(1)},\ldots,P_{\alpha((m-2)/2)}) =  \varphi_{\beta}^m(P^{'}_{\alpha(1)},\ldots,P^{'}_{\alpha((m-2)/2)})$ for all $\beta\in I(m/2, m)$.
Now we choose $\beta_0=(P_{\alpha(1)},\ldots,P_{\alpha((m-2)/2)},P_i)\in I(m/2, m)$ for $P_i\in \Sigma_m$ such that $|{\rm supp}\;\beta_0|=m$,  then  we have
 $$1=\varphi_{\beta_0}^m(P_{\alpha(1)},\ldots,P_{\alpha((m-2)/2)})=\varphi_{\beta_0}^m(P^{'}_{\alpha(1)},\ldots,P^{'}_{\alpha((m-2)/2)})$$ with

{\footnotesize $\beta_0\in\{(P^{'}_{\alpha(1)},\ldots,P^{'}_{\alpha((m-2)/2)},P_j) : |\text{supp}(P^{'}_{\alpha(1)},\ldots,P^{'}_{\alpha((m-2)/2)},P_j )| = m\},$} this implies that $\beta_0=(P^{'}_{\alpha(1)},\ldots,P^{'}_{\alpha((m-2)/2)}, P_j)$ for some $P_j \in \Sigma_m$ that is 

$P_{\alpha(1)}=P^{'}_{\alpha(1)},\ldots, P_{\alpha((m-2)/2)}=P^{'}_{\alpha((m-2)/2)}$ showing the result. 
\end{proof}

\begin{definition}
Let $m\geq 2$ be  even integer, $\overline{T}_m$  a triangulation of $C_{\frac{m-2}{2}}(\Sigma_m)$  and  $k=\frac{m+2}{2}$,  we define the matrix
$${\EuScript L}_k=\varphi^m (\overline{T}_m)$$
of order   $|\overline{T}_m |\times C^m_{m/2}$ and we call {\it the incidence matrix associated to triangulation}  
$\overline{T}_m$.
 \end{definition}
Note that the matrix ${\EuScript L}_k$ has $k$ ones in each row.\\

 Take $\displaystyle\alpha=(\alpha(1), \alpha(2),\ldots, \alpha((m-6)/2 ))\in I(m-6/2,m-2)$  arbitrary index   
     $$ T_1^1 = (P_{\alpha_{(1)}}, \dots, P_{\alpha(\frac{m-8}{2})}, P_{\alpha(\frac{m-6}{2})} ) \times C_2(\Sigma_{[\alpha(\frac{m-6}{2}),m]})$$
    $$T_{1}^{1}=P_{\alpha}\times
    C_2(\Sigma_{[\alpha(\frac{m-6}{2}),m]})=(P_{{\alpha}}\times
        R_{\frac{m-4}{2}}) \cup (P_{{\alpha}}\times
        R_{\frac{m-2}{2}}) \cup \cdots \cup (P_{{\alpha}}\times
        R_{m-1})$$
    where  $C_2(\Sigma_{[\frac{m-6}{2},m]})= R_{\frac{m-4}{2}} \cup R_{\frac{m-2}{2}} \cup\cdots \cup
    R_{m-1}$, with  
    \begin{align*}
     R_{\frac{m-4}{2}}&=\big\{ \textstyle  \big(P_{\alpha(\frac{m-4}{2})}, P_{\alpha(\frac{m-2}{2})}\big),    \big(P_{\alpha(\frac{m-4}{2})}, P_{\alpha(\frac{m}{2})}\big), \ldots,    \big(P_{\alpha(\frac{m-4}{2})}, P_{\alpha(m)}\big)    \big\},\\
 R_{\frac{m-2}{2}}&=\big\{ \textstyle  \big(P_{\alpha(\frac{m-2}{2})}, P_{\alpha(\frac{m}{2})}\big),    \big(P_{\alpha(\frac{m-2}{2})}, P_{\alpha(\frac{m+2}{2})}\big), \ldots,    \big(P_{\alpha(\frac{m-2}{2})}, P_{\alpha(m)}\big)    \big\},\\
  &\; \; \vdots \\
    R_{m-1}&=\big\{ \textstyle  \big(P_{\alpha((m-1))}, P_{\alpha (m)}\big) \big\}.
           \end{align*}
     Then the  first row of the matrix
    $\varphi^m(T_{1}^1)$ is
$$\varphi^m(P_{\alpha}, P_{\alpha(\frac{m-4}{2})}, P_{\alpha(\frac{m-2}{2})})=(\overbrace{1,\ldots,1}^k,0,\ldots,0),
$$
 the first two rows of the matrix $\varphi^m(T_{1}^1)$ are
    $$\begin{pmatrix}
    \varphi^m(P_{\alpha}, P_{\alpha(\frac{m-4}{2})}, P_{\alpha(\frac{m-2}{2})})  \\
    \varphi^m(P_{\alpha}, P_{\alpha(\frac{m-4}{2})}, P_{\alpha(\frac{m}{2})})
    \end{pmatrix}=\begin{pmatrix}
    \overbrace{1,1,\ldots,1}^k,& 0,\ldots,0,& 0,\ldots,0 \\
    1,0,\ldots,0,&\underbrace{1,\ldots,1}_{k-1},&0,\ldots,0
    \end{pmatrix},
$$
the first three rows of the matrix $\varphi^m(T_{1}^1)$  are
   $$
    \begin{pmatrix}
    \varphi^m(P_{\alpha}, P_{\alpha(\frac{m-4}{2})}, P_{\alpha(\frac{m-2}{2})})  \\
    \varphi^m(P_{\alpha}, P_{\alpha(\frac{m-4}{2})}, P_{\alpha(\frac{m}{2})})\\
    \varphi^m(P_{\alpha}, P_{\alpha(\frac{m-4}{2})}, P_{\alpha(\frac{m+2}{2})})
    \end{pmatrix}=\begin{pmatrix}
    \overbrace{1,1,\ldots,1}^k,& 0,\ldots,0,& 0,\ldots,0, & 0,\ldots,0 \\
    1,0,\ldots,0,&\underbrace{1,\ldots,1}_{k-1},&0,\ldots,0, & 0,\ldots,0\\
    0,1,\ldots,0,&1,0,\ldots,0&\underbrace{1,\ldots,1}_{k-2},& 0,\ldots,0
    \end{pmatrix},
    $$
    etcetera. As a consequence we obtain  $\varphi^m(P_{\alpha}\times R_{\frac{m-4}{2}})=A_k^0$.
    
    Now, let $I^k$ be the matrix given by the first $k$ rows of the
    identity matrix $I_{C_2^{k+1}}$,  where $\sqcup$ means
joining together side-by-side and aligning the bottoms of the
corresponding identity matrices. It is also easy to see that   
  $\varphi^m(P_{\alpha}\times R_{\frac{m-2}{2}})= I^k \sqcup A^0_{k-1}$, consequently we have the following
     \begin{align*} \varphi^m\big((P_{\alpha}\times
            R_{\frac{m-4}{2}})\cup ( P_{\alpha}\times
            R_{\frac{m-2}{2}})\big)& =A_k^0(I^k)\sqcup A_{k-1}^0,\\      
        \varphi^m\big((P_{\alpha}\times R_{\frac{m-4}{2}})\cup
       ( P_{\alpha}\times R_{\frac{m-2}{2}})\cup (P_{\alpha}\times
            R_{\frac{m}{2}})\big)&=A_k^0(I^{k+(k-1)}) \sqcup \\
        &\sqcup A_{k-1}^0(I^{k-1})\sqcup A_{k-2}^0,\\
      \vdots\quad\quad \quad\quad\quad\quad\quad & \quad\quad\quad \quad\quad\quad\vdots\\ 
        \varphi^m\big((P_{\alpha}\times R_{\frac{m-4}{2}})\cup \cdots \cup
        (P_{\alpha}\times R_{m-2})\big)&
        =A_k^0(I^{k+(k-1)+\cdots+2})\sqcup\\
       &\sqcup A_{k-1}^0(I^{(k-1)+\cdots+2}) \sqcup\\
        &\qquad
        \sqcup\cdots\sqcup A_{2}^0(I^2)\sqcup A_1^0,\\
        \varphi^m\big((P_{\alpha} \times R_{\frac{m-2}{2}})\cup \cdots \cup
       ( P_{\alpha}\times R_{m-1})\big)&
        =A_k^0(I^{k+(k-1)+\cdots+1})\sqcup\\
       & A_{k-1}^0(I^{(k-1)+\cdots+1})\sqcup\\
            & \sqcup\cdots\sqcup A_{2}^0(I^3)\sqcup
        A_1^0(I^1),
        \end{align*} 
        where for the last equality we just notice that
    $I^{k+(k-1)+\cdots+1}= I_{C_2^{k+1}}$ and similarly for the other
    $I^t$ in that formula. 
                  
   From the above,  the block stepped matrix
    \begin{gather}\label{imaL} 
\varphi^m(T_1^1) = A_k^0(I_{C^{k+1}_2})\sqcup A_{k-1}^0(I_{C^{k}_2})\sqcup\cdots\sqcup A_2^0(I_{C^{3}_2})\sqcup A_1^0(I_{C^{2}_2})
\end{gather}
is of order 
 $$\big(1+k+C_2^{k+1}\big)\times \big(C_2^{k+1}+C_2^{k}+\cdots+C_2^{3}+C_2^{2}\big)=
\big(1+k+C_2^{k+1}\big)\times C_3^{k+2}$$
 where as before $\sqcup$ means
joining together side-by-side and aligning the bottoms of the
corresponding identity matrices and filling the non-marked spaces on
the upper right blocks with zeroes. With the previous development we have

\begin{lemma}\label{LemmA_k^1}
Let $m\geq 8$ and   $\alpha=(\alpha_{(1)},\alpha_{(2)},\ldots, \alpha_{((m-6)/2)})$ 
an arbitrary index  of the  set  $I(\frac{m-6}{2}, m)$. Then $\varphi^m(T_{1}^1) = A_k^1$.
\end{lemma}

\begin{proof}
In (\ref{imaL}) we show $\varphi^m(T_1^1) = A_k^0(I_{C^{k+1}_2})\sqcup A_{k-1}^0(I_{C^{k}_2})\sqcup\cdots\sqcup A_2^0(I_{C^{3}_2})\sqcup A_1^0(I_{C^{2}_2})$ and by definition \ref{Ass-Matrix} we obtain that $\varphi^m(T_1^1)=A_k^1$.
\end{proof}

\begin{theorem}\label{TheorA_k^{k-3}}
Let $m\geq 8$ be and  $\alpha$ an arbitrary element of the index set  $I(\frac{m-6}{2}, m)$. 
Then $\varphi^m(T_i^1) = A_k^i$ for $i=1, \ldots, \frac{m+6}{2}$ and $\varphi^m(T_m) = A_k^{k-3}$.
\end{theorem}

\begin{proof}
Using the fact that $T_1^1=T_{\alpha}$ and Lemma \ref{LemmA_k^1} we have that
        \begin{align*}
        \varphi^m(T_{1}^1)&=A_k^0(I_{C_2^{k+1}})\sqcup A_{k-1}^0(I_{C_2^{k}})\sqcup \cdots\sqcup A_2^0(I_{C_2^{3}})\sqcup A_1^0(I_{C_2^{2}})= A_k^1,\\
        \varphi^m(T_{2}^1)&=
        A_k^1(I_{C_3^{k+2}})\sqcup A_{k-1}^1(I_{C_3^{k+1}})\sqcup \cdots\sqcup A_2^1(I_{C_3^{4}})\sqcup A_1^1(I_{C_3^{3}})= A_k^2,\\
        \varphi^m(T_{3}^1)&=
        A_k^2(I_{C_4^{k+3}})\sqcup A_{k-1}^2(I_{C_4^{k+2}})\sqcup \cdots\sqcup A_2^2(I_{C_4^{5}}) \sqcup A_1^2(I_{C_4^{4}})= A_k^3,\\
     \vdots   &\quad\quad\quad\quad\quad\quad\quad\quad\quad\quad\vdots\\
        \varphi^m(T_{\frac{m+6}{2}}^{1})&=
        A_k^{k-5}(I_{C_{(m-4)/2}^{m-2}})\sqcup A_{k-1}^{k-5}(I_{C_{(m-4)/2}^{m-3}})\sqcup \cdots\sqcup 
         A_1^{k-5}(I_{C_{(m-4)/2}^{(m-2)/2}})= A_k^{k-4}.
 \end{align*}      
In the same way we have $\varphi^m(T^1_{\frac{m+6}{2}}) = A_{k}^{k-4}$. Note that $$\varphi^m(T^2_{\frac{m+6}{2}}) =  
I^{C^{m-2}_{\frac{m-4}{2}}} \sqcup A_{k-2}^{k-4}.$$ 
Since
 \begin{align*}
  \varphi^m(T^1_{\frac{m+6}{2}} \cup T^2_{\frac{m+6}{2}}) & = A_k^{k-4}(I^{C_{\frac{m-4}{2}}^{m-2}}) \sqcup A_{k-1}^{k-4} \\
  \varphi^m(\cup_{i=1}^3 T^i_{\frac{m+6}{2}} ) & = A_k^{k-4}(I^{C_{\frac{m-4}{2}}^{m-2} + C_{\frac{m-6}{2}}^{m-3}}) 
  \sqcup   A_{k-1}^{k-4}(I^{C_{\frac{m-6}{2}}^{m-3}}) \sqcup A_{k-2}^{k-4} \\
    \varphi^m(\cup_{i=1}^4 T^i_{\frac{m+6}{2}}) & = A_k^{k-4}(I^{C_{\frac{m-4}{2}}^{m-2} + C_{\frac{m-4}{2}}^{m-3} + C_{\frac{m-4}{2}}^{m-4}} )   \sqcup A_{k-1}^{k-4}(I^{C_{\frac{m-4}{2}}^{m-3}} + I^{C_{\frac{m-4}{2}}^{m-4}}) \\
& \sqcup A_{k-2}^{k-4}(I^{C_{\frac{m-4}{2}}^{m-4}}) \sqcup A_{k-3}^{k-4} \\
  \vdots \quad \quad \quad \quad& \quad\quad\quad\quad\quad \quad\quad\quad \vdots \\
   \varphi^m(\cup_{i=1}^{\frac{m+4}{2}}T^i_{\frac{m+6}{2}} ) & = A_k^{k-4}(I^{C_{\frac{m-4}{2}}^{m-2} + C_{\frac{m-4}{2}}^{m-3} + \cdots + C_{\frac{m-4}{2}}^{\frac{m-4}{2}} })  \sqcup \cdots \sqcup 
  A_{1}^{k-4}(I^{ C_{\frac{m-4}{2}}^{\frac{m-4}{2}}})  \\
   & = A_k^{k-4}(I_{C_{(m-2)/2}^{m-1}})\sqcup 
    \sqcup A_{k-1}^{k-4}(I_{C_{(m-2)/2}^{m-2}})\sqcup \cdots\sqcup  A_1^{k-4}(I_{C_{(m-2)/2}^{(m-2)/2}}).  
  \end{align*}
Therefore
        \begin{align*}      
  \varphi^m\big(T_{m}\big)&
        =\varphi^m\big(T_{\frac{m+6}{2}}^{1}\cup T_{\frac{m+6}{2}}^{2}\cup \cdots \cup T_{\frac{m+6}{2}}^{\frac{m+2}{2}} \cup T_{\frac{m+6}{2}}^{\frac{m+4}{2}}\big)\\
        &=A_k^{k-4}(I_{C_{(m-2)/2}^{m-1}})\sqcup A_{k-1}^{k-4}(I_{C_{(m-2)/2}^{m-2}})\sqcup \cdots\sqcup \\
        & \sqcup A_2^{k-4}(I_{C_{(m-2)/2}^{m/2}})\sqcup  A_1^{k-4}(I_{C_{(m-2)/2}^{(m-2)/2}}),\\
        &=:A_k^{k-3}.
        \end{align*}    
\end{proof}

\begin{theorem} \label{2.3}
If   $m \geq 2$ is  an even integer and $k=\frac{m+2}{2}$,
then  $C_{\frac{m-2}{2}}(\Sigma_m)=T_m$ is a \it{triangulated set}  and  ${\EuScript L}_k=  A_k^{k-3}$. 
 \end{theorem} 

\begin{proof}
First we show that $C_{\frac{m-2}{2}}(\Sigma_m)$ is a triangulated set. If $m=2,4,6$ is obvious. Now we suppose that $m \geq 8$, then 
clearly $\overline{T}_m\subseteq C_{\frac{m-2}{2}}(\Sigma_m)$, for all $m$ even integer. Moreover  $\varphi^m$ is an injective function 
see Lemma \ref{Funt-Inj}, and by Theorem \ref{TheorA_k^{k-3}}, $\varphi^m\big(T_{m}\big)=A_k^{k-3}$ so then we have  $\vert T_m\vert=\vert (Im\; \varphi^m\vert _{T_m})\vert $, equal to the number of rows in the matrix $A_k^{k-3}$  and we obtain
$$|T_m|=C^{k+(k-1)-1}_{(k-1)-1}=C^{2k-2}_{k-2}=C^m_{\frac{m-2}{2}}$$ that it implies $C_{(m-2)/2}(\Sigma_m)=T_m$,  and varying $(P_{\alpha(1)},\ldots,P_{\alpha((m-2)/2)})\in T_m$  we have the matrix
    $${\EuScript L}_k=\varphi^m(T_m)=  A_k^{k-3}. $$
\end{proof}

\begin{corollary} If $k$ is any positive integer, then 
$$A_k^{k-2} = {\EuScript L}_k(I_{C_{k-1}^{2k-2}})\sqcup \;^t{\EuScript L}_k,$$
$$ {\EuScript L}_k  = A_k^{k-4}(I_{C_{k-2}^{2k-3}})\sqcup A_{k-1}^{k-3}.$$
\end{corollary}

\begin{proof}  From Theorem \ref{exis} and \ref{trans}
\begin{align*}
A_k^{k-2}&  = A_k^{k-3}(I_{C_{k-1}^{2k-2}})\sqcup A_{k-1}^{k-2} \\
& = {\EuScript L}_k(I_{C_{k-1}^{2k-2}})\sqcup \; ^t A_{k}^{k-3}\\
& = {\EuScript L}_k(I_{C_{k-1}^{2k-2}})\sqcup \;^t{\EuScript L}_k.
\end{align*} 
The second equality is followed in a similar way.
\end{proof}

\begin{theorem}
Suppose that the coefficients of the matrix $A_k^{k-2}$ is in ${\mathbb{F}}$, an arbitrary field of 
char$({\mathbb{F}}) = 0$ or char$({\mathbb{F}})\geq r$ where $2 \leq k \leq r$ for $r= \lfloor \frac{n+2}{2} \rfloor$. Then
the matrix $A_k^{k-2}$ has maximum range equal to $C_{k-1}^{2k-1}$.
\end{theorem}

\begin{proof}
We have $A_k^{k-2}$ is square and fragmented, i. e., 
$$ A_k^{k-2} =  {\EuScript L}_k(I_{C_{k-1}^{2k-2}})\sqcup \;^t{\EuScript L}_k.$$
Clearly the submatrix $I_{C_{k-1}^{2k-2}} \sqcup \;^t{\EuScript L}_k$ has maximum range.  Also the submatrix  $ {\EuScript L}_k$, block in  $A_k^{k-2}$, has maximum rank as shown in Theorem 9 of \cite{2}.
In addition, neither row of submatrix $I_{C_{k-1}^{2k-2}} \sqcup \;^t{\EuScript L}_k$ can be placed as a linear combination of rows of submatrix (block) ${\EuScript L}_k $ and reciprocally. Therefore the maximum range of $A_k^{k-2}$ is $C_{k-1}^{2k-1}$.
\end{proof}

\section{The matrix $\mathcal{M}$ as a combinatorial design}\label{Secc5}
In this section we show the combinatorial construction of $\mathcal{M}$, other sparse $(0,1)$-matrix using the matrix  $A_k^{k-3}$.
The study is divided into two cases even and  odd.

\subsubsection{Case: $n$ even}
 For $n\geq 4$ even and $r:=\frac{n+2}{2}$ consider integers $1\leq k \leq r-2$ and sequences of integers
 $1\leq a_1<a_2<\cdots<a_{2k}\leq 2n$ such that $a_i+a_j\neq 2n+1$. 
Defining
   $$\Sigma_{a_1,\ldots,a_{2k}}:=\{P_i\in\Sigma_n: i+a_j\neq 2n+1,\, 2n-i+1+a_j\neq 2n+1\}$$
we have that: 
    \begin{enumerate}
        \item $\big|\Sigma_{a_1,\ldots,a_{2k}}\big|=n-2k$.
        \item If $k \geq 1$, setting
      {\footnotesize  
       $\Sigma\{a_1,\ldots,a_{2k}\}:=(a_1,\ldots, a_{2k})\times C_{(n-2(k+1))/2}(\Sigma_{a_1,\ldots ,a_{2k}}) \\
        =\big\{\big(a_1,\ldots, a_{2k}, P_{\alpha(1)},\ldots, P_{\alpha((n-2(k+1))/2)}\big):\\
      \big(\alpha(1),\ldots, \alpha((n-2(k+1))/2)\big) \in I(
        (n-2(k+1))/2,n-2k)\big\},$ 
        }
        where $1 \leq k \leq r-2$.       
   \end{enumerate}

With the above notation we have.

\begin{lemma}\label{Set-Trian}
If $n\geq 4$ is even, $r=\frac{n+2}{2}$ and $k=1,\ldots, r-2$,  then the set 
 
 $C_{(n-2(k +1))/2}(\Sigma_{a_1,\ldots ,a_{2k}})$  is a triangulated set.
\end{lemma}

\begin{proof}
We get to take $m=n-2k$ in the definition \ref{TrianG}. 
\end{proof}

\begin{lemma}\label{lem2.5}
If $n\geq 4$ is even,  then
 \begin{equation*}
I(n-2,2n)=C_{(n-2)/2}(\Sigma_n)\cup\Big(\bigcup_{k=1}^{r-2}
\bigcup_{\substack{(a_1,\ldots, a_{2k})\in I(2k, 2n)\\ a_i+a_j\neq 2n+1}}\Sigma\{a_1,\ldots, a_{2k}\}\Big).
\end{equation*}
is a partition of the set $I(n-2,2n)$.
\end{lemma}

\begin{proof}
It is sufficient to show that 
 $$I(n-2,2n)\subset C_{(n-2)/2}(\Sigma_n)\cup\Big(\bigcup_{k=1}^{r-2}
\bigcup_{\substack{(a_1,\ldots, a_{2k})\in I(2k, 2n)\\ a_i+a_j\neq 2n+1}}\Sigma\{a_1,\ldots, a_{2k}\}\Big).$$
Let $\overline{a}=(a_1,\ldots, a_{n-2})\in I(n-2, 2n)$ such that $a_i+a_j\neq 2n+1$
for all $i, j$, then $\overline{a}\in \Sigma\{a_1,\ldots, a_{n-2}\}.$
On the other hand, if $\overline{a}=(a_1,\ldots, a_{n-2})=(P_{\alpha(1)},\ldots, P_{\alpha((n-2)/2)})$ up a permutations in its components, and where $P_{\alpha(i)}=(a_{\alpha(i)}, a_{\alpha(2n-i+1)})$  so this case   $\overline{a}\in C_{(n-2)/2}(\Sigma_n)$. 
Finally for $\overline{a}=(a_1, \ldots, a_{2k}, P_{\alpha(1)},\ldots, P_{\alpha(\frac{n-2(k+1)}{2})})$ where $a_i+a_j\neq 2n+1$ and $\alpha(i)+a_j\neq 2n+1$ then
$\overline{a}\in \Sigma\{a_1,\ldots, a_{2k}\}$. With this we show the desired.
\end{proof}

\subsubsection{Case: $n$ odd}

 The odd case is obtained by modifications to the even case.

   Explicitly,  for $n\geq 5$ odd and $r=(n+1)/2$, consider integers $1\leq k \leq
r-2$ and sequences of integers
 $1\leq a_1<a_2<\cdots<a_{2k +1}\leq 2n$ such that $a_i+a_j\neq 2n+1$, and define
    $\Sigma_{a_1,\ldots,a_{2k +1}}=\{P_i\in\Sigma_n: i+a_j\neq 2n+1,\, 2n-i+1+a_j\neq 2n+1\}$.
Then:\medskip

 \begin{enumerate}
\item  $\big|\Sigma_{a_1,\ldots,a_{2k + 1}}\big|=n-(2k +1)$.
\item  If $k \geq 1$, setting
   {tiny     
        $\Sigma\{a_1,\ldots,a_{2k + 1}\}=(a_1,\ldots, a_{2k + 1})\times C_{(n-(2k+3))/2}(\Sigma_{a_1,\ldots ,a_{2k +1}}) \\
        =\big\{\big(a_1\ldots a_{2k +1}, P_{\alpha(1)},\ldots, P_{\alpha((n-(2k+3))/2)}\big)$
          }
such that $\big(\alpha(1),\ldots, \alpha((n-(2k+3))/2)\big)\in I(
        (n-(2k+3))/2,n-(2k +1))\big\}$ for $1\leq k
        \leq r-2$. 
\item Let $\Sigma(i)=\{P_j\in \Sigma_n: i\neq j \}$.
\end{enumerate}       
      
\begin{lemma}
If $n\geq 5$ is odd,
 $r=\frac{n+1}{2}$ and $k=1,\ldots, r-2$,  then the set 

  $C_{(n-(2k+3))/2}(\Sigma_{a_1,\ldots ,a_{2k +1}})$ is a triangulated set.
\end{lemma}

\begin{proof}
We get to take $m=n-(2k+3)$ in the definition \ref{TrianG}.  
\end{proof}

\begin{lemma}\label{lem3.5}
Let $n\geq 5$ odd with the above notation we have
{\scriptsize\begin{equation}\label{lem2.51}
        I(n-2,2n)= \left(\bigcup_{i=1}^n \left[(i)\times C_{\frac{n-3}{2}}(\Sigma(i))
        \right]\right)  \cup \Bigg(\bigcup_{k=0}^{r-2}\;\; \bigcup_{\substack{1\leq
            a_1<\ldots< a_{2k +1} \leq 2n \\ a_i+a_j\neq 2n+1}}
        \Sigma\{a_1, \ldots, a_{2k+1}\} \Bigg).
\end{equation} }
is a partition of the set $I(n-2,2n)$.
\end{lemma}
\begin{proof}
The demonstration is similar to Lemma  \ref{lem2.5}.
\end{proof}

\subsection{An incidence function over  the index set $I(n-2, 2n)$}

First, we extend the function $\varphi^m$ defined in \eqref{Sec-Mm} to $I(n-2, 2n)$.\\
 For $n\geq 4$,  we consider
$$\varphi:I(n-2,2n)\rightarrow \{0,1\}^{I(n-2, 2n)} $$ 
such that  {\footnotesize $(\alpha(1),\ldots,\alpha(n-2))\mapsto \Big( \varphi_{\beta}(\alpha(1),\ldots,\alpha(n-2))\Big)_{\beta\in I(n,2n)}$}\\
where {\scriptsize
$\varphi_{\beta}(\alpha(1),\ldots,\alpha(n-2)) = \begin{cases}
1  & \text{if}\; \beta\in\{(i,\alpha_{rs}, 2n-i+1) \\  & : |\text{supp} \{(i,\alpha_{rs}, 2n-i+1) \}| = n\}, \\
0& \text{otherwise}.
\end{cases}$} 
 \begin{definition}\label{Mdesing}
 Let $n\geq 4$ be, the matrix
\begin{gather*}\label{MatrixM}
    {\mathcal{M}}=\Big(\varphi_{\beta}({\alpha_{rs}})\Big)_{\alpha_{rs}\in I(n-2,2n),\beta\in I(n,2n)}.
    \end{gather*}
of order $C^{2n}_{n-2}\times C^{2n}_n$ will be called \it{the incidence matrix}.
 \end{definition}   

\begin{lemma}\label{lemmaTri}
Let $n\geq 4$ be. Then $\varphi$ is injective, also 
 \begin{center} 
$\varphi|_{C_{(n-2)/2}(\Sigma_n)}=\varphi^n$ and $\varphi|_{ \Big((i)\times C_{\frac{n-3}{2}}(\Sigma(i))\Big)}=\varphi^n$.
\end{center} 
\end{lemma}    
\begin{proof}
It is enough to show  the injectivity.  

So suppose $\varphi\big(\alpha(1),\ldots,\alpha(n-2))= \varphi\big(\alpha(1)^{'},\ldots,\alpha(n-2)^{'}\big)$, then 
$ \varphi_{\beta}(\alpha(1),\ldots,\alpha(n-2))= \varphi_{\beta}(\alpha(1)^{'},\ldots,\alpha(n-2)^{'})$ for all $\beta\in I(n, 2n)$. If $\beta_0=(i, \alpha(1),\ldots,\alpha(n-2), 2n-i+1)$ for some $1\leq i \leq 2n$ such that $|\text{supp} \;\beta_0|=n$, then 
$$1= \varphi_{\beta_0}(\alpha(1),\ldots,\alpha(n-2))= \varphi_{\beta_0}(\alpha(1)^{'},\ldots,\alpha(n-2)^{'})$$  implies that 
$\beta_0=(j, \alpha^{'}(1),\ldots,\alpha^{'}(n-2), 2n-j+1)$ and so $\alpha(1)=\alpha^{'}(1),\ldots, \alpha(n-2)=\alpha^{'}(n-2)$ which concludes our demonstration.
\end{proof}
    
  \begin{lemma}\label{FuntMa}
  For $n\geq 4$ positive integers,  let $r= \frac{n+2}{2}$ be if $n$ is even or let  $r= \frac{n+1}{2}$  if $n$ odd.  
   Then
  \begin{enumerate}
\item $\varphi(C_{(n-2)/2}(\Sigma_n))={\EuScript L}_r$.
\item $\varphi(T_{n-2k}^{a_1,\ldots ,a_{2k}})={\EuScript L}_{r-k}^{a_1,\ldots ,a_{2k}}$ where  
${\EuScript L}_{r-k}^{a_1,\ldots ,a_{2k}}={\EuScript L}_{r-k}$ with $k=1,\ldots, r-2$ .
\item $ \varphi \Big((i)\times C_{\frac{n-3}{2}}(\Sigma(i))\Big)={\EuScript L}^i_r$, a copy of ${\EuScript L}_r$ for all $i=1,\ldots, n$.
\item $\varphi(T_{n-(2k+1)}^{a_1,\ldots ,a_{2k+1}})={\EuScript L}_{r-k}^{a_1,\ldots ,a_{2k+1}}$ where ${\EuScript L}_{r-k}^{a_1,\ldots ,a_{2k+1}}={\EuScript L}_{r-k}$with $k=0,\ldots, r-2$.
\end{enumerate}
  \end{lemma}  
 \begin{proof}
  If $n$ is even and $k=0,\ldots, r-2$, then   
  $C_{\frac{(n-2k)-2}{2}}(\Sigma_{a_1,\ldots ,a_{2k}})=T_{n-2k}$ is 
  triangulated set, by Theorem \ref{2.3}  $\varphi(T_{n-2k})=A_{\frac{(n-2k)+2}{2}}^{\frac{(n-2k)+2}{2}-3}={\EuScript L}_{\frac{n+2}{2}-k}={\EuScript L}_{r-k}$. 
In the same way, if $n$ is odd and  $k=0,\ldots, r-2$, then $C_{\frac{(n-(2k+1))-2}{2}}(\Sigma_{a_1,\ldots ,a_{2k +1}})=T_{n-(2k+1)}$ is triangulated set, by Theorem \ref{2.3} $\varphi(T_{n-(2k+1)})=A_{\frac{n-(2k+1)+2}{2}}^{\frac{n-(2k+1)+2}{2}-3}={\EuScript L}_{\frac{n-(2k+1)+2}{2}}={\EuScript L}_{r-k}$. 
 \end{proof}
    
\begin{theorem} \label{thmPeven}
    Let $n\geq 4$ be an even integer,  $r=(n+2)/2$ and $1\leq k\leq r-2$.  Then $\mathcal{M}$ is  a  direct sum, that is
    $${\mathcal{M}}={\EuScript L}_r\oplus\Bigg( \bigoplus_{k=1}^{r-2}\Big(\bigoplus_{\substack{1\leq a_1<\cdots<a_{2k}\leq 2n\\ a_i+a_j\neq 2n+1}}{\EuScript L}_{r-k}^{(a_1,\cdots, a_{2k})}\Big)\Bigg),$$
    where ${\EuScript L}_{r-k}^{(a_1,\cdots, a_{2k})}$ is a copy of ${\EuScript
        L}_{r-k}$, for each $1 \leq k \leq r-2$.
\end{theorem} 

\begin{proof}    
    By Lemma \ref{lem2.5} 

\begin{align*}
        I(n-2,2n)&=
        C_{(n-2)/2}(\Sigma_n)\cup\Bigg(\bigcup_{k=1}^{r-2}\bigcup_{\substack{(a_1,\ldots,
            a_{2k})\in I(2k,2n)\\ a_i+a_j\neq 2n+1}}\Sigma\{a_1,\ldots, a_{2k}\}\Bigg),
\end{align*}
using the definition \ref{Mdesing} and (1) of Lemma \ref{FuntMa}  we have
{\tiny \begin{align*}
        {\mathcal{M}} &= (\varphi(C_{\frac{n-2}{2}}(\Sigma_n))) \bigoplus
        \Bigg(  \bigoplus_{k-1}^{r-2} \bigoplus_{\substack{(a_1,\ldots,a_{2k})\in I(2k,2n) \\ a_i+a_j\neq 2n+1}} (\varphi(\Sigma(a_1,\ldots,\alpha_{2k})))\Bigg)\\
        {\mathcal{M}} & = {\EuScript L}_r\oplus\Bigg( \bigoplus_{k=1}^{r-2}\Big(\bigoplus_{\substack{1\leq a_1<\cdots<a_{2k}\leq 2n\\ a_i+a_j\neq 2n+1}}\!\!\!\!{\EuScript L}_{r-k}^{(a_1,\cdots, a_{2k})}\Big)\Bigg).
\end{align*}} 
\end{proof}

\begin{theorem} \label{thmPodd}
    Let $n\geq 5$ be an odd integer  and $r=(n+1)/2$. Then  $\mathcal{M}$  is a    direct sum, that is
     $${\mathcal{M}}={\EuScript L}_r^{n}\oplus\Bigg(\bigoplus_{k=1}^{r-2}
    \Bigg(\bigoplus_{\substack{1\leq a_1< a_2< \cdots < a_{2k+1}\leq 2n \\ a_i
        + a_j \neq 2n+1}} {\EuScript
        L}_{r-k}^{(a_1,a_2,\dots,a_{2k+1})}\Bigg)\Bigg),$$ 
        where ${\EuScript L}_{r-k}^{(a_1, a_2,\dots, a_{2k+1})}$ is a copy of ${\EuScript L}_{r-k}$
    for each $1\leq k\leq r-2$ and ${\EuScript L}_r^n = {\EuScript L}_r \oplus \cdots \oplus {\EuScript L}_r$ $n$-times.
\end{theorem} 

\begin{proof}
    As before, let $\Sigma_n=\{P_1,\dots,P_n\}$ and
    $\Sigma(i)=\Sigma_s-\{P_i\}$ for all $i \in \{1,\ldots, n\}$ and
    $1\leq s \leq n$. Then, from (2) of  Lemma \ref{FuntMa}, we have that it is easy
    to see that the image of $(i)\times
    C_{\frac{n-3}{2}}\big(\Sigma(i)\big)$ under $\varphi$ is ${\EuScript
        L}_r$, i.e., 
  \begin{equation}
 \varphi \Big((i)\times
    C_{\frac{n-3}{2}}\big(\Sigma(i)\big)\Big)={\EuScript L}_r^{i}\;\;\text{\rm a copy of}\;\; 
    {\EuScript L}_r \;\text{for\; all}\; i=1,\dots, n.
 \end{equation}    
    Similarly, by  definitions of $\varphi$ and   $\Sigma\{a_1,\ldots,a_{2\ell + 1}\}$ we have that
 \begin{equation}\label{eq12}
 \varphi(\Sigma( a_1, a_2,\ldots, a_{2k+1}))= {\EuScript  L}_{r-k}^{(a_1, a_2,\ldots, a_{2k+1})}.
  \end{equation}     
Now,  using  Lemma \ref{lem2.51}, we obtain that  
    $$    I(n-2,2n)=   \left( \bigcup_{i=1}^{n}\left[ (i)\times
        C_{\frac{n-3}{2}}\big(\Sigma(i)\big)\right] \right)   \cup
        \Bigg(
        \bigcup_{k=1}^{r-2}\;\bigcup_{\substack{1\leq a_1  < \dots < a_{2k+1}
            \leq 2n \\ a_i + a_j \neq
            2n+1}}\Sigma\{a_1,\dots,a_{2k+1}\} \Bigg),
$$
 where $\Sigma\{a_1,\dots,a_{2k+1}\}=(a_1,a_2,\dots,a_{2k+1})\times
    C_{\frac{n-(2k+3)}{2}}\big(\Sigma_{a_1,a_2,\dots,a_{2k+1}}\big)$.
    Using   definition \ref{Mdesing}   and  (2) of Lemma \ref{FuntMa}, we obtain
\begin{align*} {\mathcal{M}}& =\varphi(I(n-2,2n))\\
      & =  \bigoplus_{i=1}^n \varphi \Big((i)\times C_{\frac{n-3}{2}}(\Sigma
        (i))\Big) \oplus  \bigoplus_{k=1}^{r-2}\Bigg ( \bigoplus_{\substack{1\leq a_1<
            a_2< \ldots < a_{2k+1}\leq 2n \\ a_i + a_j \neq 2n +1}}
        \varphi(\Sigma( a_1, a_2,\ldots, a_{2k+1}))\Bigg)   \\
        &=\bigoplus_{i=1}^n {\EuScript L}_r^i \oplus
        \Bigg(\bigoplus_{k=1}^{r-2} \Bigg(\bigoplus_{\substack{1\leq a_1 < \ldots <
            a_{2k+1} \leq a_{2n} \\ a_i + a_j \neq 2n+1}}  {\EuScript
            L}_{r-k}^{(a_1, a_2,\ldots, a_{2k+1})}\Bigg)\Bigg).
 \end{align*}
\end{proof}

\section{APPLICATIONS}\label{Secc6}
\subsection{Canonical matrix of the Lagrangian-Grassmannian variety}

Throughout this subsection, let  $E$ be a $2n$--dimensional vector space over $\mathbb{F}$, an arbitrary field, 
equipped with a linear symplectic structure $\langle\;,\;\rangle$ see \cite{0.01}. 
Then the Lagrangian-Grassmannian variety is
$$L(n, E)= Z\big(Q_{\alpha^{'}, \beta^{'}}, \Pi_{\alpha_{st}}\big),$$ 
where
  $Z(\cdot)\subseteq {\mathbb P}\big(\bigwedge^{\ell}\overline{E}\big)$ denotes of zeros the Pl\"ucker relations $Q_{\alpha^{'}, \beta^{'}}$ see \ref{eq1.1} in Appendix A  and of the linear homogeneous polynomials $\Pi_{\alpha_{st}}$, here  $\alpha^{'} \in I(n-1, m)$, $\beta^{'} \in I(n+1,m)$ respectively $\alpha_{st}\in I(n-2,2n)$. 
Now  we define {\it the support} of $\alpha=(\alpha_1,\ldots,\alpha_{\ell}) \in I(\ell,m)$  as
the set $\text{supp}\;\alpha=\{\alpha_1,\ldots,\alpha_{\ell}\}$ and for all $\alpha_{st}\in I(n-2,2n)$ we define the
a linear polynomial $\Pi_{\alpha_{st}}\in {\mathbb F}[X_{\alpha}:\alpha\in I(n,2n)]$, as 
\begin{gather*}\label{eq2.1}
\Pi_{\alpha_{st}}:=\sum_{i=1}^nc_{i,\alpha_{st},2n-i+1}X_{i,\alpha_{st},2n-i+1}
\end{gather*}
with
$$c_{i,\alpha_{st},2n-i+1}=\begin{cases}
1  & \text{if $|\text{supp}\{i,\alpha_{st},2n-i+1\}|=n$}, \\
0 & \text{otherwise}.
\end{cases}$$

Let $B_{L(n,E)}$  be the matrix  of order  $C^{2n}_{n-2}\times C^{2n}_{n}$  associated the system of homogeneous linear equations  $$\Pi:=\{\Pi_{\alpha_{st}}: \alpha_{st}\in I(n-2,2n)\}.$$

\begin{lemma}
Let  $E$ symplectic vector space of dimension $2n$ and  $\mathcal{M}$ the matrix defined in the Section \ref{Secc5}. Then $B_{L(n, E)}={\mathcal{M}}$.
\end{lemma}
\begin{proof}
By definition of matrix $B_{L(n, E)}$ each row is of the form $(c_{\beta})_{\beta\in I(n,2n)}$  where there is a $\alpha_{rs}\in I(n-2, 2n)$such that 
$$c_{\beta} = \begin{cases}
1  & \text{if}\; \beta\in\{(i,\alpha_{st}, 2n-i+1)  : |\text{supp} \{(i,\alpha_{st}, 2n-i+1) \}| = n\}, \\
0& \text{otherwise}.
\end{cases}$$ \\
this implies that for every row of $B_{L(n, E)}$ we have  $(c_{\beta})_{\beta \in I(n,2n)}=\varphi({\alpha_{st}})$ for some 
$\alpha_{st}\in I(n-2, 2n)$, 
therefore $B_{L(n, E)}={\mathcal{M}}$ see definition \ref{Mdesing}.
\end{proof}

\begin{lemma}\label{sse}
Keeping the above notation. Let ${\mathbb{F}}$ be a field such that char$({\mathbb{F}}) = 0$ or
char$({\mathbb{F}}) \geq r$ where $r= \lfloor \frac{n+2}{2} \rfloor$.
\begin{enumerate}
\item[A)] If $H$ is a matrix of order $C_{n-2}^{2n} \times C_n^{2n}$ and maximum rank that annuls the rational points of 
$L(n,E)(\mathbb{F}_q)$, then $H= PB_{L(n,E)}$, here $P$ is an invertible matrix.

\item[B)] Suppose that there exists $R$ matrix such that $L(n,E) = G(n,E) \cap \ker R$. Then $R=PB_{L(n,E)}$ where $P$ is an invertible matrix.
\end{enumerate}    
 \end{lemma}
\begin{proof}
By Theorem 6 of \cite{2} we have  $rank B_{L(n, E)}=C^n_{n-2}$ and therefore $\langle \Pi_{\alpha_{rs}}: \alpha_{rs}\in I(n-2, 2n) \rangle_{{\mathbb F}}$ is a vector space of dimension $C^n_{n-2}$. Now let  $H= ( h_1, h_2, \ldots,  h_{\epsilon})^t$ be the matrix of rank $\epsilon=C^{2n}_{n-2}$ where $\{h_1,\ldots, h_{\epsilon}\}\subset (\bigwedge^nE)^*$ are its rows. As  
$H$  annuls the rational points of $L(n, E)({\mathbb F}_q)$, by the Theorem 12 of \cite{2} we have $\langle h_1, \dots, h_{\epsilon}\rangle_{{\mathbb F}} \subset\langle \Pi_{\alpha_{rs}}: \alpha_{rs}\in I(n-2, 2n) \rangle_{{\mathbb F}}$ and  $h_1,\ldots, h_{\epsilon}$ also forms a basis,  therefore $H = PB_{L(n, E)}$ where $P$ is the non-singular matrix of basis change of the symplectic vector space $E$.
 Now suppose that  $R=(h_1, h_2, \ldots h_{\epsilon})^t$
 is a rank matrix $\epsilon$  such that $L(n, E)=G(n, E)\cap \ker R$, then $L(n, E)({\mathbb F}_q)\subset \ker R$ and $\epsilon \leq C^n_{n-2}$. If $\epsilon=C^n_{n-2}$ the affirmation is followed by what was said before.
Now suppose that $t<C^n_{n-2}$ then $\ker B_{L(n, E)} \varsubsetneq \ker R$ this implies that
$$
L(n, E) =G(n, E)\cap \ker B_{L(n, E)} 
           \varsubsetneq G(n, E) \cap \ker R 
           = L(n, E)$$
 which is a contradiction and therefore $\epsilon =C^n_{n-2}$.
\end{proof}

\begin{theorem}\label{MatrixB} 
Let ${\mathbb F}$ a field such  that char$({\mathbb {F}})=0$ or  char$({\mathbb F})\geq r$ where $r=\lfloor  \frac{n+2}{2}\rfloor$. Then the matrix $B_{L(n, E)}$ it is diagonal by blocks and is unique, up basis change, that annuls the rational points of the Lagrangian-Grassmannian variety.
\end{theorem}
\begin{proof}
 From Theorems \ref{thmPeven} and \ref{thmPodd} the matrix $B_{L(n, E)}={\mathcal{M}}$  is a direct sum of the matrices ${\EuScript L}_k$ with $k=2,\ldots, r$ and $r=\lfloor  \frac{n+2}{2}\rfloor$, even more it is easy to see that diagonal by blocks. By  Lemma  \ref{sse} we have the rest of this theorem.
\end{proof}

\subsection{Construction of  Isodual Codes}
Another class of sparse linear codes is the low density generator matrix (LDGM) codes, which have sparse generator matrices. In this section we build two classes of isodual codes, one class is from LDGM-codes and the other class is `dense" codes, we will show that they share properties. Both codes are defined, based on their code generating matrix, which is obtained as following method, which we summarize for reasons of space:\\

\begin{description}
\item[{\rm Step 1. }]  Using the algorithm \ref{Algoritm1} 
(respectively the algorithm \ref{Algoritm 2}) build the spread matrix $A_{k-1}^{k-3}$, (respectively the dense matrix $B_{k-1}^{k-3}$)
\item[{\rm Step 2.} ] For the matrix $A_{k-1}^{k-3}$we "add" to the left the matrix  $I_{C_{k-2}^{2k-3}}$ to generate  $\big( I_{C_{k-2}^{2k-3}}| A_{k-1}^{k-3}\big)$(respectively we add to the left the matrix $J_{C_{k-2}^{2k-3}}$ to generate the dense matrix $\big(J_{C_{k-2}^{2k-3}}| B_{k-1}^{k-3}\big)$). 
\end{description}

\begin{remark}
Note that these matrices are obtained by truncating the matrices. ${\EuScript  L}_{k}$ (respectively $ {\mathcal {M}}_k$).
\end{remark}

We denote the matrices resulting from the algorithms above by:
$${\EuScript  L}_{k}^{truc}:=\big( I_{C_{k-2}^{2k-3}}| A_{k-1}^{k-3}\big)\;\;and\;\; {\mathcal{M}}_{k}^{truc}:=\big(J_{C_{k-2}^{2k-3}}| B_{k-1}^{k-3}\big).$$
We also define the matrices 
$$H_s:=\big(\big( A_{k-1}^{k-3} \big)^t | I_{C_{k-2}^{2k-3}}\big)\;\;and \;\; H_d =\big(\big( B_{k-1}^{k-3} \big)^t | I_{C_{k-2}^{2k-3}}\big),$$
which are the check parity matrices  ${\EuScript  L}_{k}^{truc}$  and ${\mathcal{M}}_{k}^{truc}$, respectively.
We denote by $C_{{\EuScript  L}_{k}^{truc}}$ and  $C_{{\mathcal {M}}_{k}^{truc}}$  the  $[2\dot C_{k-2}^{2k-3}, \dot C_{k-2}^{2k-3}]$ codes generated by matrices ${\EuScript  L}_{k}^{truc}$ and ${\mathcal{M}}_{k}^{truc}$ respectively. Analogously we denote by $C_{H_s}$ and  $C_{H_d}$ the codes generated by the matrices $H_s$ and $H_d$.
 Clearly $C_{H_s}=C_{{\EuScript  L}_{k}^{truc}}^{\perp}$ and $C_{H_d}=C_{{\mathcal{M}}_{k}^{truc}}^{\perp}$
\begin{lemma}\label{IsodualCode}
The codes $C_{{\EuScript  L}_{k}^{truc}}$ and $C_{{\mathcal{M}}_{k}^{truc}}$  are   isoduales .
\end{lemma}
\begin{proof}
By  Corollary \ref{cuaPe}
we have $\big( A_{k-1}^{k-3} \big)^t=I^a_{C_{k-2}^{2k-3}}A_{k-1}^{k-3} I^a_{C_{k-2}^{2k-3}}$.
 Analogously by  Corollary \ref{cuaPe2} we have 
 $\big( B_{k-1}^{k-3} \big)^t=I^a_{C_{k-2}^{2k-3}}B_{k-1}^{k-3} I^a_{C_{k-2}^{2k-3}}$. Then
\begin{align*}
H_s & =\big(\big( A_{k-1}^{k-3} \big)^t | I_{C_{k-2}^{2k-3}}\big)\\
& = \big(I^a_{C_{k-2}^{2k-3}}A_{k-1}^{k-3} I^a_{C_{k-2}^{2k-3}}| \big( I_{C_{k-2}^{2k-3}}\big)\big)\\
& = I^a_{C_{k-2}^{2k-3}} \big( I_{C_{k-2}^{2k-3}}| A_{k}^{k-3} \big)I^a_{C_{k-2}^{2k-3}}\\
&\sim \big( I_{C_{k-2}^{2k-3}}| A_{k}^{k-3}\big)={\EuScript  L}_{k}^{truc}
\end{align*}
Which implies that $C_{{\EuScript  L}_{k}^{truc}}\sim C_{{\EuScript  L}_{k}^{truc}}^{\perp}$, so
$C_{{\EuScript  L}_{k}^{truc}}$ is isodual.\\
The proof that  $C_{{\mathcal{M}}_{k}^{truc}}$ is isodual, is similar.
\end{proof}
\begin{theorem}
Let $C_{{\EuScript  L}_{k}^{truc}}$ and $C_{{\EuScript  M}_{k}^{truc}}$ the $[2C_{k-2}^{2k-3}, C_{k-2}^{2k-3}]$-binary codes and  $k$  a positive integer even.  Then
$$
W_{C_{{\EuScript  L}_{k}^{truc}}}(x,y) = W_{C_{{\EuScript  M}_{k}^{truc}}}(x,y) 
=\sum_{i=0}^{\lfloor \frac{1}{4}C_{k-2}^{2k-3} \rfloor}a_ig_1(x,y)^{C_{k-2}^{2k-3}-4i}g_2(x,y)^i,
$$
where $g_1(x, y)=y^2+x^2$ and $g_2(x, y)=y^8+14x^4y^4+x^8$ are Gleason polynomials, see \cite{3.5}.
Moreover
$$d_{C_{{\EuScript  L}_{k}^{truc}}}=d_{C_{{\EuScript  M}_{k}^{truc}}}\leq\begin{cases}
2 \lfloor   \frac{n_0}{4}\rfloor +2  & if \;\; n_0\leq 30, \\
2 \lfloor   \frac{n_0}{4}\rfloor & if \;\; n_0\geq 32, \\
\end{cases}$$
here $n_0:=C_{k-2}^{2k-3}$.
\end{theorem}
\begin{proof}
From Lemma \ref{IsodualCode} and Theorem 9.2.1 (i) in \cite{3.5},  it is enough to see that the codes are even.
It is easy to see that the generating matrix 
${\EuScript  L}_{k}^{truc}:=\big( I_{C_{k-2}^{2k-3}}| A_{k-1}^{k-3}\big)$  respectively   ${\mathcal {M}}_{k}^{truc}:=\big(J_{C_{k-2}^{2k-3}}| B_{k-1}^{k-3}\big)$ they have lines of weight $k$, respectively $C_{k-2}^{2k-3}-k$, which are positive pairs because by hypothesis $k$ is even, and by the Theorem 1.4.11 in \cite{3.5} the codes 
 $C_{{\EuScript  L}_{k}^{truc}}$ and $C_{{\mathcal{M}}_{k}^{truc}}$ they are of even weight.
\end{proof}
\subsection{Efficient encoders based on approximate lower diagonal matrix}
Given $H$ a sparse matrix, the parity-check code $C_H$ associated with $H$ is called {\it Low Density Parity Check} (LDPC) code and is a linear error correcting code see \cite{3}.
In this subsection we shall develop an  encoder for  $C_{A_k^{\ell-2}}$ LDPC-code regular generated by   $A_k^{\ell-2}$ sparce,  fragmented and is in approximate lower diagonal matrix. This encoder, is efficient because it is based on `Efficient encoders based on approximate lower triangulations" developed by Richardson-Urbanke  see \cite{4}.

Let $k$ and $\ell$ be any positive integers, by Theorem  \ref{exis}, the matrix $A_k^{\ell-2} = A_k^{\ell-3}(I_{C_{\ell -1}^{k+\ell -2}}) \sqcup A_{k-1}^{\ell-2}$ is of order   $C^{k+\ell-1}_{\ell-1}\times C^{k+\ell-1}_{\ell}$, and by definition of $A_k^{\ell-2}$ we have  
$$
 A_k^{\ell-2} =A_k^{\ell-3}\left( I_{C^{k+\ell-2}_{\ell-1}}
\right)\sqcup A_{(k-1)}^{\ell-3}\left( I_{C^{k+\ell-3}_{\ell-1}}
\right) \sqcup \cdots \sqcup  \sqcup A_2^{\ell-3}\left(
I_{C^{\ell}_{\ell-1}} \right) \sqcup A_1^{\ell-3}\left(
I_{C^{\ell-1}_{\ell-1}} \right) .$$

We wrote
$$A_k^{\ell-2}=\begin{bmatrix}
A_k^{\ell-3}\;\; & 0\;\; & 0 \;\; \\
 I_{C^{k+\ell-2}_{\ell-1}}\;\;& B\;\; &  A \;\;
\end{bmatrix}$$
where $B\sqcup A=A_{(k-1)}^{\ell-3}\left( I_{C^{k+\ell-3}_{\ell-1}}
\right) \sqcup \cdots \sqcup A_2^{\ell-3}\left(
I_{C^{\ell}_{\ell-1}} \right) \sqcup A_1^{\ell-3}\left(
I_{C^{\ell-1}_{\ell-1}} \right)$,  with
$B$  of order $C^{k+\ell-2}_{\ell-1}\times C^{k+\ell-2}_{\ell-2}$ and  $A$  of order  $C^{k+\ell-2}_{\ell-1}\times (C^{k+\ell-1}_{\ell}-C^{k+\ell-1}_{\ell-1})$.

We consider a solution of  type $x=(P_2, P_1, S)$,
 $P_2$  is an order vector $C^{k+\ell-2}_{\ell-1}$, $P_1$  is an order vector 
$C^{k+\ell-2}_{\ell-2}$, and $S$  is an order vector $(C^{k+\ell-1}_{\ell}-C^{k+\ell-1}_{\ell-1})$. As in \cite{4}, $g:=C^{k+\ell-2}_{\ell-2}$ denote the {\it gap}, and measures in some way to be made precise shortly, the  ''distance''  of the given parity-check matrix to a lower triangular matrix. 

Multiplying this matrix by the left 
 $$ \begin{bmatrix} 
0\;\; &\;\;  I_{C^{k+\ell-2}_{\ell-1}}  \\
I_{C^{k+\ell-2}_{\ell-2}} & -A_k^{\ell-3}  
\end{bmatrix}
\begin{bmatrix}
A_k^{\ell-3}\;\; & 0\;\; & 0 \;\; \\
 I_{C^{k+\ell-2}_{\ell-1}}\;\;& B\;\; &  A \;\;
\end{bmatrix}
 =
\begin{bmatrix} 
-BP_1^t-AS^t\;\; & B\;\;\;  & \;\;\;  A  \\
0 &-A_k^{\ell-3} B &\;\;\;-A_k^{\ell-3} A 
\end{bmatrix}$$
which multiplied by $x^t=(P_2, P_1, S)^t$, where $\Phi^t$ define transpose matrix,
gives us the vector $$\big( P_2+BP_1^t+AS^t,\;\;\; -A_k^{\ell-3}BP_1^t-A_k^{\ell-3}AS^t \big)=(0,0).$$ 
Let $\phi=-A_k^{\ell-3}B$ the matrix of order $C^{k+\ell-2}_{\ell-2}\times C^{k+\ell-2}_{\ell-2}$ and suppose  that is invertible, we can simply perform further column permutations to remove this singularity, see section II of \cite{4}. 
Then $P_2=-BP_1^t-AS^t$ and  \;$P_1=\phi^{-1}(A_k^{\ell-3}A)S^t$.

\begin{theorem}
   With the previous notation. Let $A_k^{\ell-2}$, sparce,  fragmented and  approximate lower diagonal matrix
and  $S\in {\mathbb F}^{\{C_{\ell}^{k+\ell-1}-C_{\ell-1}^{k+\ell-1}\}}$  an arbitrary vector.
Then  $S\mapsto (P_2, P_1, S)$ is an encoder for $C_{A_k^{\ell-2}}$  LDPC-code regular. 
\end{theorem}
\begin{proof}
Remember that $P_1=\phi^{-1}(A_k^{\ell-3}A)S^t$ and   $P_2=-BP_1^t-AS^t$ then multiplying by the right we have
 {\scriptsize \begin{align*}
  & \begin{bmatrix} 
A_k^{\ell-3}\;\; & 0\;\;\;  & \;\;\; 0  \\
 I_{C^{k+\ell-2}_{\ell-1}} & B &\;\;\; A 
\end{bmatrix}
\begin{bmatrix}
- BP_1^t-AS^t \\
 \phi^{-1}(A_k^{\ell-3}A)S^t\\
 S^t
\end{bmatrix}
 =  \\
 & =   \begin{bmatrix} 
-A_k^{\ell-3}BP_1^t-A_k^{\ell-3}AS^t, &  \big( - BP_1^t-AS^t \big)+\big(B \phi^{-1}(A_k^{\ell-3}A)S^t \big) + AS^t
\end{bmatrix}\\ 
& = \begin{bmatrix} 
\phi P_1^t-\phi P_1^t, & -BP_1^t+ BP_1^t
\end{bmatrix} =[0,0],
\end{align*}} 
thus the theorem is proved.
\end{proof}


\section{Conclusion}
 We present a link between three areas of mathematics; algorithms, linear code theory and the geometry of 
 the Lagrangian-Grassmannian variety, where the key ingredient are the LG-matrices.
 The key to all of this is an algorithm that builds
  structured sparse $(0, 1)$-matrix in approximate lower triangular form, 
with $k$ ones in each row and $\ell$ ones in each column, where $k$ and $\ell$ are any positive integers.
We study some properties of the two families of $(0,1)$-matrices built with our algorithm.
We give a new construction of  isodual lineal codes and dense codes, and we study both codes.
The $ (0,1) $-matrices presented in this work can have various applications.

\section*{ Appendice: Lagrangian--Grassmannian variety }\label{append-A}
For a vector space $E$ of finite dimension $m$ over arbitrary field $F$, let
$G(\ell,m)$ denote the Grassmannian variety of vector subspaces of dimension $\ell$ of $E$. The
{\it Pl\"ucker embedding}  $\rho: G(l,m) \rightarrow {\mathbb P}(\wedge^{\ell} E)$ maps a subspace $W\in G(\ell,m)$ to $\rho(W)=v_1\wedge \cdots \wedge v_{\ell}$, where $\{v_1,\ldots, v_{\ell}\}$ is a basis of $W$. If 
$\{e_1,\ldots,e_m\}$ is a basis of $E$, and for $\alpha=(\alpha_1,\ldots,
\alpha_{\ell})\in I(\ell,m)$  we put $e_\alpha = e_{\alpha_1}
\wedge \cdots \wedge e_{\alpha_{\ell}}$,  then the set
$\{ e_\alpha\}_{\alpha\in I(\ell,m)}$ is a basis of $\wedge^{\ell} E$. Now, writing 
 $w \in \wedge^{\ell} E$ as $w= \sum_{\alpha
\in I(\ell,m)} P_\alpha e_\alpha$, the scalars $P_\alpha$ are called  the {\it
Pl\"ucker coordinates} of $w$ and $w_\rho:=(P_\alpha)_{\alpha \in I(n,2n)}$ is
 the  {\it Pl\"ucker vector} of $w$. 
Now, if $w = \sum_{\alpha \in I(\ell,m)}P_\alpha e_\alpha \in
{\mathbb{P}}(\wedge^{\ell} E)$, then $w\in G(\ell,m)$ if and only if for
each pair of tuples $1\leq \alpha_1< \cdots <
\alpha_{\ell-1}\leq m$ and $1\leq \beta_1< \cdots < \beta_{\ell+1}\leq m$, the Pl\"ucker coordinates of $w$ satisfy the quadratic called Pl\"ucker relations see \cite{6}
\begin{equation}\label{eq1.1}
 Q_{\alpha^{'}, \beta^{'}}:=\sum_{i=1}^{\ell+1} (-1)^i P_{\alpha_1\cdots \alpha_{\ell-1}\beta_i}P_{\beta_1\beta_2\cdots
\widehat{\beta_i} \cdots \beta_{\ell+1}} = 0,
\end{equation}
where  $\widehat{\beta_i}$ means that the corresponding term is omitted and where $\alpha^{'} \in I(n-1, m)$,                  $\beta^{'} \in I(n+1,m)$.
A special case  is when the  vector space $E$  over the  field $F$ has a skew-symmetric nondegenerate bilinear (symplectic)  form 
$\langle \,, \,\rangle$. If $2n$ is the dimension of $E$,  there is a basis ${\mathfrak B}=\{ e_1, \ldots ,
e_{2n}\}$ of $E$, where 
$$ \langle e_i , e_j \rangle  =
\begin{cases}
  1 & \text{if},      j=2n-i+1, \\
   0   & \text{otherwise}.
\end{cases}
$$

We define the Lagrangian--Grassmannian variety 
$$L(n,2n)=\{ v_1\wedge \cdots \wedge v_n\in G(n,2n)\; :\; \langle v_i, v_j \rangle=0 \} $$ 
as a subset of the  Grassmannian $G(n,2n)$. It was shown in \cite{1} that   $L(n,2n)=G(n,2n)\cap {\mathbb P}(\ker f)$, where $f$ is the contraction map.
 For $\alpha=(\alpha_1,\ldots, \alpha_r,\ldots,\alpha_s,\ldots, \alpha_k)\in I(k,2n)$ let 
 $\alpha_{rs}=(\alpha_1,\ldots, \widehat{\alpha}_r,\ldots,\widehat{\alpha}_s,\ldots, \alpha_k)\in I(k-2,2n)$ deleting the numbers below the hat symbol. Following \cite{1}  define the linear polynomials    
\begin{equation}\label{eq2.2}
\Pi_{\alpha_{rs}}:=\sum_{i=1}^nc_{i,\alpha_{rs},2n-i+1}X_{i,\alpha_{rs},2n-i+1}
\end{equation}
with $X_{i,\alpha_{rs},2n-i+1}$  indeterminates and 
$$c_{i,\alpha_{rs},2n-i+1}=\begin{cases}
1  & \text{if $|\text{supp}\{i,\alpha_{rs},2n-i+1\}|=k$}, \\
0 & \text{otherwise},
\end{cases}$$
where $\text{supp}(\beta)=\{\beta_1,\ldots,\beta_d\}$ for
$\beta=(\beta_1,\ldots,\beta_d)\in I(d,2n)$. 
In \cite{3} it was shown that $\ker f = Z\big ( \Pi_{\alpha_{rs}}: \alpha_{rs}\in I(n-2, 2n)\big)\subseteq{\mathbb P}\big(\bigwedge^{\ell}\overline{E}\big)$.
So then we can conclude that the Lagrangiana-Grassmannian variety is the set
$$L(n, E)= Z\big(Q_{\alpha^{'}, \beta^{'}}, \Pi_{\alpha_{st}}\big)$$ 
where
  $Z(\cdot)\subseteq {\mathbb P}\big(\bigwedge^{\ell}\overline{E}\big)$ denotes the set of zeros of the given polynomials,
  $\alpha^{'} \in I(n-1, m)$, $\beta^{'} \in I(n+1,m)$, $\alpha_{st}\in I(n-2,2n)$.  
  Moreover the set of rational points is defined by 
  $$L(n, E)({\mathbb F}_q)= Z\big(Q_{\alpha^{'}, \beta^{'}}, \Pi_{\alpha_{st}}, X^q_{\alpha}-X_{\alpha}\big)_{\alpha\in I(n,2n)}$$ 
 In \cite{2}, it is shown that if $h\in  \big(\bigwedge^{\ell}\overline{E}\big)^{*}$ suct that $h(L(n, E)({\mathbb F}_q))=0$, then $h\in \langle\Pi_{\alpha_{rs}} | \alpha_{rs}\in I(n-2, 2n) \rangle$.

\section*{Acknowledgment}

The authors were supported by Laboratorio de
Cifrado y Codificaci\'on de la Informaci\'on (LCCI-UACM) of
Universidad Aut\'onoma de la Ciudad de M\'exico.

\end{document}